\documentclass[narrowdisplay,seceqn,secthm]{elsart1p}

\usepackage{amsmath,amssymb,float}
\usepackage{amsfonts}
\usepackage{graphics,euscript,epsfig}
\usepackage{enumerate}

\usepackage[english]{babel}
\usepackage[T1]{fontenc}
\usepackage[latin1]{inputenc} 
\usepackage{textcomp}
\usepackage{ae}
\usepackage{aecompl}
\usepackage[cm]{aeguill}
\usepackage{lmodern}

\newtheorem{thrm}{Theorem}[section]

\newtheorem{coro}[thrm]{Corollary}

\newtheorem{rmrk}[thrm]{Remark}

\newcommand{\R}{\ensuremath{\mathbb{R}}}

\newcommand{\N}{\ensuremath{\mathbb{N}}}

\numberwithin{equation}{section}
\numberwithin{figure}{subsection}
\everymath{\displaystyle}

\begin{document}

\begin{frontmatter}



\title{On the 
    Korteweg-de Vries approximation\\ for uneven bottoms\vspace{0.7em}}


\author{Florent Chazel}

\address{\vspace{0.1em}\\Saint-Venant Laboratory for Hydraulics\\(Université Paris-Est, Joint research unit EDF R\&D - CETMEF - Ecole des Ponts)\\6 quai Watier, BP 49, F-78401 Chatou, France}

\ead{florent-externe.chazel@edf.fr}

\begin{abstract}
In this paper we focus on the water
waves problem for uneven bottoms on a
two-dimensionnal domain. Starting from the symmetric Boussinesq
systems derived in [Chazel, {\it Influence of bottom topography on long water waves}, 2007], we recover the uncoupled
Korteweg-de Vries (KdV) approximation justified in
[Schneider-Wayne, {\it The long-wave limit for the water-wave problem. I. The case of zero surface tension}, 2002] for flat bottoms, and in [Iguchi, {\it A long wave approximation for capillary-gravity waves and an effect of the bottom}, 2007] in the context of
bottoms tending to zero at infinity at a substantial rate. The goal of this paper is to investigate the validity of this approximation for more general bathymetries. We exhibit two kinds of topography for which this approximation
diverges from the Boussinesq solutions. A
topographically modified KdV approximation is then proposed to deal
with such bathymetries, where topography-dependent terms are added to the solutions of the KdV equations. Finally, all the models involved are numerically computed and compared.
\end{abstract}

\end{frontmatter}

\section*{Introduction}

The water waves problem consists in describing the evolution of the
free surface and velocity field of a layer of ideal, incompressible
and irrotationnal fluid under the only influence of gravity. The governing
equations - also called free surface Euler equations - are fully 
non-linear and non-strictly hyperbolic, and their direct study and 
computation remains a real obstacle. Many
authors such as Nalimov (\cite{Nalimov}, 1974), Yoshihara
(\cite{Yoshihara1}, 1982), Craig (\cite{Craig}, 1985), Wu (\cite{Wu1},
1997 and \cite{Wu2}, 1999), Ambrose-Masmoudi (\cite{AM}, 2005) and Lannes (\cite{Lannes}, 2005 and 
\cite{AL1}, 2007) have
successfully tackled the problem of well-posedness of these equations.
Nevertheless, the numerical computation of these solutions
remains a tough task, especially in 3-D - see the works of Grilli et al. (\cite{Grilli}, 2001) and Fochesato-Dias (\cite{Foch}, 2001). \\
An alternative way to describe these solutions and their time
behaviour is to look for approximations via the use of asymptotic models. Such models are usually derived formally from the water waves problem by introducing dimensionless parameters. 
Making some hypothesis on these parameters reduces the framework to
more limited physical regimes but allows the construction of
asymptotic models. In this work, we focus on the so-called long waves
regime. In this regime, the ratios $\varepsilon=a/h_0$ and
$\mu=h_0^2/\lambda^2$ where $a$ denotes the typical amplitude of the
waves, $h_0$ the mean depth and $\lambda$ the typical wavelength, are
small and of the same order. Many models can be found in the
litterature corresponding to this regime. Among them, we can quote the
works of Boussinesq (\cite{Boussinesq1}, 1871 and \cite{Boussinesq2},
1872) who was the first to propose a model that take into account both
nonlinear and dispersive effects, the unidirectional models such as
the Korteweg-de Vries (KdV) equation (\cite{KdV}, 1895), the
Kadomtsev-Petviashvili (KP) equation (\cite{KP}, 1970) and the
Benjamin-Bona-Mahony one (\cite{BBM}, 1972). These historical models
have been considerably studied and generalized, and their justification has been investigated among others by Craig (\cite{Craig}, 1985), Schneider-Wayne (\cite{SW}, 2000), Bona-Colin-Lannes (\cite{BCL}, 2005), Lannes-Saut (\cite{LannesJCS}, 2006) for flat bottoms, and Iguchi (\cite{Iguchi}, 2006), Chazel (\cite{FC}, 2007), Alvarez-Lannes (\cite{AL1}, 2007) for uneven bottoms.\\
If we focus more specifically on the KdV approximation, we emphasize that the investigation of this model on uneven bottoms is not recent : this subject has been tackled among others by Ostrovskii-Pelinovskii (\cite{Ostrovskii}, 1970), Kakutani (\cite{Kakutani}, 1971), Johnson (\cite{Johnson}, 1973), Miles (\cite{Miles}, 1979) and Newell (\cite{Newell}, 1985).  However, their approach is quite different compared to the one chosen here, as will be clarified later in this paper. The more recent articles of Iguchi (\cite{Iguchi}, 2006) and Chazel (\cite{FC}, 2007) are somehow the starting point of the present work. In \cite{Iguchi}, Iguchi derived two different models : a coupled KdV system for relatively general bottom topographies, and a uncoupled KdV system for bathymetries decaying at a substantial rate at infinity. In \cite{FC}, the author derived two classes of symmetric Boussinesq systems for two bottom topography scales, corresponding respectively to slightly and largely varying bottoms. The aim of this work is to recover the uncoupled KdV model - also justified by Schneider-Wayne \cite{SW} for flat bottoms - starting from any previous Boussinesq system proposed in \cite{FC} for slightly varying bottoms, to discuss its validity regarding the bottom topography, and to propose an alternative.

\subsection*{Formulation of the problem}

In this paper, we work in two dimensions : $x$ corresponds to
the horizontal coordinate and $y$ to the vertical one. We
denote by $(t,x) \rightarrow \eta(t,x)$ and $x \rightarrow
b(x)$ the parametrizations of the free surface and bottom, defined respectively 
over the surface $y=0$ and the mean depth $y=-h_0$ at the steady state. 
The time-dependant domain $\Omega_t$ of the fluid is thus taken of the form : 
$$\Omega_t = \{(x,y),\,x\in\R,\,-h_0+b(x) \le y\le\eta(t,x)\}\;\;.$$
\numberwithin{figure}{section}
\begin{figure}[t]
\caption{Representation of the fluid domain}
\center{\epsfig{file=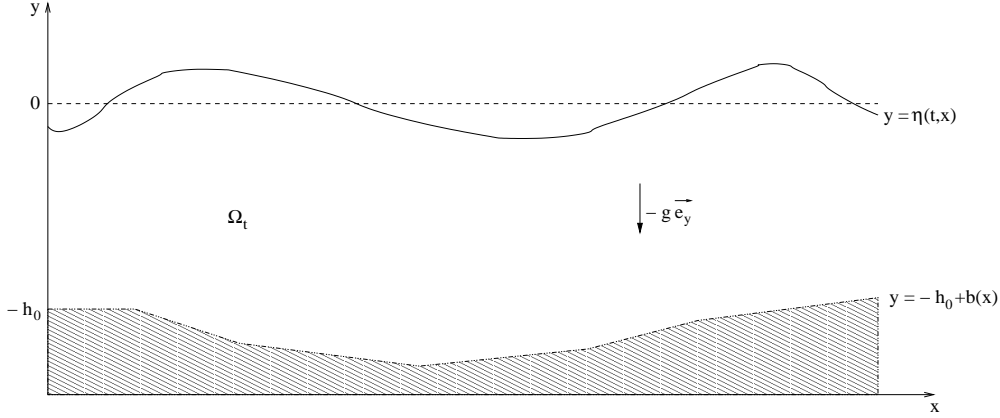,width=13.5cm,height=5.5cm}}
\end{figure}
\numberwithin{figure}{subsection}
\noindent For the sake of simplicity, we assume here that $b \in W^{k,\infty}(\R^d)$, $k$ being as    	  large as  needed, where we recall that $W^{k,\infty}(\R^d) = \{u\,/\,\partial_x^m u \in L^{\infty}(\R^d), 0 \le m \le k\}$. In order to avoid some special physical cases
such as the presence of islands or beaches, we set a condition of minimal water depth : there exists
a strictly positive constant $h_{min}$ such that
\begin{equation} \label{hmin}
\eta(t,x) + h_0 - b(x)  \ge h_{min}\;\;,\;\; (t,x) \in \R\times\R\;\;.
\end{equation}
We introduce the following dimensionless parameters :
$$
\varepsilon = \frac{a}{h_0}\;\;,\;\;\mu=\frac{h_0^2}{\lambda^2}\;\;,\beta=\frac{\beta_0}{h_0}\;\;,
$$
where $a$ denotes the typical amplitude of the waves, $h_0$ the mean depth, $\lambda$ the typical wavelength and $\beta_0$ the typical amplitude of the bottom topography. In the present long waves regime, one has $\varepsilon \ll 1$, $\mu \ll 1$ and $\varepsilon \approx \mu$, i.e. the Stokes number $S=\varepsilon/\mu$ is of order $O(1)$. Moreover, one has $\beta = B_0 \varepsilon$ with $B_0 = O(1)$ since we focus here on slightly varying bottoms, i.e. bathymetries of small amplitude. We do not make any mild-slope assumption on the spatial variations of the bottom, i.e. the bathymetry can vary rapidly. For the sake of simplicity, we take the Stokes number $S$ and the constant $B_0$ equal to one, which implies that we have $\mu = \beta = \varepsilon$. This choice only lightens the writings by suppressing the constants $S$ and $B_0$ from the equations, and has no influence on the following results. 

\vspace{1em}

In \cite{FC}, the author justified a whole class of
symmetric Boussinesq systems as being asymptotic models to the water
waves problem for slightly varying bottoms in 2-D and 3-D, assuming
the existence of the water waves solutions on a large time scale. This
assumption has been recently proved by Alvarez-Lannes in \cite{AL1},
where the authors systematically justified the main asymptotics models
used in coastal oceanography. The justification of these symmetric
Boussinesq models is hence complete. All the details on the construction and
justification of these models can be found in \cite{FC}; 
the expression of such symmetric systems in 1-D surface is as follows :
$$
(\Sigma)
\left\{
\begin{array}{l}
\vspace{0.5em}
(1-\varepsilon a_2 \partial_x^2) \partial_t v + \partial_x
\eta + \varepsilon \Big[ \frac{1}{2} \eta \partial_x \eta +
  \frac{3}{2} v \partial_x v - \frac{1}{2} b \partial_x \eta + a_1
  \partial_x^3 \eta \Big] = 0\,, \\
(1-\varepsilon a_4 \partial_x^2) \partial_t \eta +
\partial_x v + \varepsilon \Big[ \frac{1}{2} \partial_x \left(
    (\eta-b) v \right) + a_3 \partial_x^3 v \Big] = 0\,,
\end{array}
\right.
$$
with 
$$
\left\{
\begin{array}{lll}
\vspace{0.2em}
a_1 = -\lambda_1\frac{\theta^2-1}{2} & ; &
a_2 = (\lambda_1-1)\frac{\theta^2-1}{2} \,,\\
a_3 = \lambda_2 (\frac{\theta^2}{2}-\frac{1}{6}) & ; &
a_4 = (1-\lambda_2) (\frac{\theta^2}{2}-\frac{1}{6}) \,,
\end{array}
\right.
$$
and $(\theta,\lambda_1,\lambda_2)$ being chosen in $[0,1]\times\R^2$ such that
$a_1=a_3$, $a_2\ge 0$, $a_4 \ge 0$. The parameter $\theta$ determines the height $y = -1 + \theta (\varepsilon\eta +1 - \beta b)$ where the above velocity $v$ is taken. The parameters $\lambda_1$ and $\lambda_2$ are arbitrary parameters and have no physical meaning. \\
Several different set of values for $(\theta,\lambda,\mu)$ can be
found such that the previous condition on $a_1,a_2,a_3,a_4$ is
verified, f.e. $\theta = \sqrt{2/3}$, $\lambda_1=\lambda_2=1/2$ which gives $a_1=a_2=a_3=a_4=1/12$. Throughout this paper, we denote by $(\Sigma)$ any of this
symmetric system. This system is the starting point of our work.
	
\subsection*{Organization of the paper}

The paper is organized as follows. In Section I, we recover the KdV approximation proposed by Schneider-Wayne \cite{SW} and Iguchi \cite{Iguchi} as follows : we first diagonalize the system $(\Sigma)$ and then look for an approximation of the solution of the obtained coupled system. This approximation is searched under the form of a couple of waves $(U_0,N_0)$ moving in opposite directions, plus correcting terms $(U_1,N_1)$ that satisfy a sublinear growth condition. We show that $U_0$ and $N_0$ must satisfy a system of two uncoupled KdV equations, as shown by Iguchi in \cite{Iguchi}. It is proved that this approximation is correct for sufficiently decaying initial data and topography.\\
In Section II, we discuss the validity of this approximation - which
is equivalent to check the sublinear growth condition on the
correctors $(U_1,N_1)$ - for non trivial bathymetries. An analysis of
these terms allows us to weaken the decay assumptions on the
initial data and the bathymetry. However, it appears clearly that some
more general bathymetries can unvalidate the model, and two examples
of bottoms are provided for which the approximation diverges : a
bottom corresponding to a simple step and a slowly varying
sinusoidal bottom. A topographically modified KdV
approximation is then proposed by adding bottom-dependent correcting
terms to the classical approximation.\\
In Section III, all these models - the Boussinesq one, the usual
uncoupled KdV one and the topographically modified version - are
numerically integrated and compared on the two bathymetries introduced
in Section II. The numerical schemes are presented, and the results show that both the Boussinesq model and the alternative KdV approximation successfully reproduce the expected physical phenomenons.

\section{The classical KdV approximation}

In this section, we recover the usual uncoupled KdV approximation
justified by Schneider and Wayne in \cite{SW} for flat botoms, and
by Iguchi in \cite{Iguchi} for slightly varying bottoms. To this
end, we start with any of the symmetric Boussinesq system ($\Sigma$) 
derived in \cite{FC} and look for approximate solutions of the diagonalized 
version of ($\Sigma$) under the form of two waves moving in opposite
directions. Each wave is shown to be slightly modulated in time with a dynamic
governed by a KdV equation, while the correcting terms must
solve an inhomogeneous transport equation. The sum of these two waves is proved
to give an approximation of the solutions of ($\Sigma$) for sufficiently decaying initial
data and topography.

\subsection{Derivation of the approximation}

Let us recall the expression of ($\Sigma$) :
$$
(\Sigma)
\left\{
\begin{array}{l}
\vspace{0.5em}
(1-\varepsilon a_2 \partial_x^2) \partial_t v + \partial_x
\eta + \varepsilon \Big[ \frac{1}{2} \eta \partial_x \eta +
  \frac{3}{2} v \partial_x v - \frac{1}{2} b \partial_x \eta + a_1
  \partial_x^3 \eta \Big] = 0\,, \\
(1-\varepsilon a_4 \partial_x^2) \partial_t \eta +
\partial_x v + \varepsilon \Big[ \frac{1}{2} \partial_x \left(
    (\eta-b) v \right) + a_3 \partial_x^3 v \Big] = 0\,,
\end{array}
\right.
$$
with $a_2 \ge 0$ and $a_4 \ge 0$. We point out that such a system is well posed for
sufficiently smooth initial data and provides an approximation of the
water waves problem of order $O(\varepsilon)$ for times of order
$O(1/\varepsilon)$ (see \cite{FC} and \cite{AL1} for further
details).\\
In order
to recover the KdV approximation, we first diagonalize $(\Sigma)$ by
introducing the following unknowns :
$$
U = v + \eta \;\;;\;\; N = v-\eta\;.
$$
Plugging the relations $v = (U+N)/2$ and $\eta = (U-N)/2$
into $(\Sigma)$ yields the following coupled system
$(\Gamma)$ in terms of $U$ and $N$ :
$$
(\Gamma)
\left\{
\begin{array}{l}
\vspace{0.5em}
\partial_t U + \partial_x U + \varepsilon \left[
  \frac{1}{8}\partial_x (3 U^2+N^2+2 U N ) - \frac{1}{2} b
  \partial_x U - \frac{1}{4} \partial_x b (U+N) \right. \\ 
\vspace{0.5em}
  \hspace{7em} \left. - \frac{a_2}{2}\partial_x^2\partial_t (U+N) - \frac{a_4}{2}
  \partial_x^2\partial_t (U-N) + a_1 \partial_x^3 U \right] = 0\,,\\
\vspace{0.5em}
\partial_t N - \partial_x N + \varepsilon \left[
  \frac{1}{8}\partial_x (U^2 + 3 N^2 + 2 U N ) + \frac{1}{2} b
  \partial_x N + \frac{1}{4} \partial_x b (U+N) \right. \\
  \hspace{7em} \left. - \frac{a_2}{2}\partial_x^2\partial_t (U+N) + \frac{a_4}{2}
  \partial_x^2\partial_t (U-N) - a_1 \partial_x^3 N \right] = 0\,.
\end{array}
\right.
$$
At this step, we choose to look for an
approximate solution $(U_{app}, N_{app})$ of $(\Gamma)$ of the form :
\begin{eqnarray} \label{ansatz}
\left\{
\begin{array}{l}
\vspace{1em}
U_{app}(t,x) = U_0(T, x-t)+\varepsilon U_1(T,t,x)\;, \\
N_{app}(t,x) = N_0(T, x+t)+\varepsilon N_1(T,t,x)\;,
\end{array}
\right.
\end{eqnarray}
where $T$ is the slow time variable $T=\varepsilon t$. The choice of variables $x-t$ for $U_0$ and $x+t$ for $N_0$ emphasizes that we look for solutions which - at first order - propagate at speed one, to the right for $U$ and to the left for $N$ : this comes from the first order equations on $U$ and $N$ which are respectively $(\partial_t+\partial_x)U = 0$ and $(\partial_t-\partial_x)N = 0$.\\
The use of two different time scales is useful in capturing both the short time evolution of the wave and the nonlinear and dispersive dynamics which occur for larger time scales.\\
We complement this ansatz with two initial conditions on $(U_{app}, N_{app})$ and a classical sublinear growth condition on the correctors $(U_{1}, N_{1})$ :
$$
U_{app\,|_{t=0}} = U_{|_{t=0}} \;\;,\;\; N_{app\,|_{t=0}} = N_{|_{t=0}}\;,
$$
and for all $T_0 \ge 0$
\begin{eqnarray}\label{growth}
\left\{
\begin{array}{l}
\vspace{1em}
\lim_{t\rightarrow +\infty} \frac{1}{t}\,\left|\,U_1(.\,,t,.)\right|_{L^{\infty}([0,T_0];L^2(\R))} = 0
\;, \\
\lim_{t\rightarrow +\infty} \frac{1}{t}\,|\,N_1(.\,,t,.)|_{L^{\infty}([0,T_0];L^2(\R))} = 0\;.
\end{array}
\right.
\end{eqnarray}
Such a condition on the correctors is quite usual in multiscales expansions and has been first introduced in the context of nonlinear geometric optics by Joly, Métivier and Rauch in \cite{JMR2}. It forces the correctors $\varepsilon U_1$ and $\varepsilon N_1$ to be small on the large time scale associated to the KdV dynamics ; namely, it ensures that $\varepsilon U_1 = o(1)$ and $\varepsilon N_1 = o(1)$ in $L^{\infty}([0,\frac{T_0}{\varepsilon}];L^2(\R))$.

\vspace{1em}

\noindent Plugging this ansatz into $(\Gamma)$ and neglecting the terms of order $O(\varepsilon)$ yields the following system :
\begin{equation} \label{coupled}
\left\{
\begin{array}{l}
\vspace{0.5em}
(\partial_t + \partial_x) U_1 = f(T,x-t) - \frac{1}{8} \partial_x N_0^2 - \frac{1}{4} \partial_x (U_0 N_0) + \frac{1}{2} b \partial_x U_0 + \frac{1}{4} U_0 \partial_x b \\
\hspace{7em} + \frac{1}{4} N_0 \partial_x b + \frac{a_2-a_4}{2}\, \partial_x^3 N_0\;, \\
(\partial_t - \partial_x) N_1 = g(T,x+t) - \frac{1}{8} \partial_x U_0^2 - \frac{1}{4} \partial_x (U_0 N_0) - \frac{1}{2} b \partial_x N_0 - \frac{1}{4} U_0 \partial_x b \\
\hspace{7em} - \frac{1}{4} N_0 \partial_x b - \frac{a_2-a_4}{2}\, \partial_x^3 U_0\;,
\end{array}
\right.
\end{equation}
\begin{equation}
\mbox {where} \;\; \left\{
\begin{array}{l}
\vspace{0.5em}
f(T,x-t) = -\partial_T U_0 - \frac{3}{8} \partial_x U_0^2 - (a_1+\frac{a_2+a_4}{2})\, 
\partial_x^3 U_0\;,\\
g(T,x+t) = -\partial_T N_0 - \frac{3}{8} \partial_x N_0^2 + (a_1+\frac{a_2+a_4}{2})\, \partial_x^3 N_0\;.
\end{array}
\right.
\end{equation}
\underline{N.B.} : For the sake of simplicity, we have kept here the notations $\partial_x$ on $U_0$ and $N_0$ while we should have written rigorously $\partial_{X_{-}}U_0$ and $\partial_{X_{+}}N_0$ with $X_{-}=x-t$ and $X_{+}=x+t$.

\vspace{1em}

\noindent At this step, an explicit resolution of (\ref{coupled}) in
terms of $U_0$ and $N_0$ shows that $U_1$ and $N_1$ have the
simplified form :
$$
\left\{
\begin{array}{l}
\vspace{0.5em}
U_1(T,t,x) = t\times f(T,x-t) + h_b^1(U_0(T,x-t),N_0(T,x+t))\;, \\
N_1(T,t,x) = t\times g(T,x+t) + h_b^2(U_0(T,x-t),N_0(T,x+t))\;,
\end{array}
\right.
$$
the complete expression of $(U_1,N_1)$ being given at the end of this
section in (\ref{explicit}).\\
Such expressions of $U_1$ and $N_1$ include terms that grow linearly
in time, which is inconsistent with the sublinear growth conditions
(\ref{growth}). It follows that both $f(T,x-t)$ and $g(T,x-t)$ must be
null quantities for all $(T,t,x) \in [0,T_0]\times\R^2$. Writing
explicitly this result, remarking that $a_1+(a_2+a_4)/2 =
(a_1+a_3+a_2+a_4)/2 = 1/6$ and plugging it into (\ref{coupled}), one obtains the following uncoupled KdV equations on $(U_0,N_0)$ and inhomogeneous transport equations on the correctors $(U_1,N_1)$ : 
\begin{equation}
(\Sigma_{KdV})\;\left\{
\begin{array}{l}
\vspace{0.5em}
\partial_T U_0 + \frac{3}{8}\, \partial_x U_0^2 + \frac{1}{6}\, \partial_x^3 U_0 = 0\;,\\
\partial_T N_0 + \frac{3}{8}\, \partial_x N_0^2 - \frac{1}{6}\, \partial_x^3 N_0 = 0\;,
\end{array}
\right.
\end{equation}
and
\begin{equation}
(\Sigma_{corr})\;\left\{
\begin{array}{l}
\vspace{0.5em}
(\partial_t + \partial_x) U_1 = - \frac{1}{8} \partial_x N_0^2 - \frac{1}{4} \partial_x (U_0 N_0) + \frac{1}{2} b \partial_x U_0 + \frac{1}{4} U_0 \partial_x b \\
\hspace{7em} + \frac{1}{4} N_0 \partial_x b + \frac{a_2-a_4}{2}\, \partial_x^3 N_0\;, \\
(\partial_t - \partial_x) N_1 = - \frac{1}{8} \partial_x U_0^2 - \frac{1}{4} \partial_x (U_0 N_0) - \frac{1}{2} b \partial_x N_0 - \frac{1}{4} U_0 \partial_x b \\
\hspace{7em} - \frac{1}{4} N_0 \partial_x b - \frac{a_2-a_4}{2}\, \partial_x^3 U_0\;.
\end{array}
\right.
\end{equation}
Finally, we can construct an approximation of the solutions of the initial system $(\Sigma)$ in a natural way : let $(v^{\,\varepsilon}_{\Sigma},\eta^{\,\varepsilon}_{\Sigma})_{0\le\varepsilon\le\varepsilon_0}$ be a family of solutions of $(\Sigma)$ with initial data $(v_0,\eta_0)$. One defines $u_0 = v_0+\eta_0$ and $n_0 = v_0-\eta_0$, and let $(U_0,N_0)$ be the solutions of the uncoupled KdV equations $(\Sigma_{KdV})$ with initial data $(u_0,n_0)$. In the end, the uncoupled KdV approximation of the solutions of $(\Sigma)$ is given by :
\begin{eqnarray} \label{M}
(\mathcal{M})\;\;:\;\;
v_{KdV}^{\;\;\;\varepsilon} = \frac{U_0+N_0}{2} \;\;;\;\;
\eta_{KdV}^{\;\;\;\varepsilon} = \frac{U_0-N_0}{2}\;.
\end{eqnarray}
This approximation corresponds to the one proposed by Schneider-Wayne \cite{SW} for flat bottoms and by Iguchi \cite{Iguchi} for uneven bottoms.

\begin{rmrk} 
Iguchi derived his KdV approximation in the framework of capillary-gravity waves, and not gravity waves as in this paper. However, we can see in \cite{Iguchi} that the only impact of this capillarity-gravity waves approach lies in the - constant - coefficient in front of the dispersive terms. The comparison between this work and the present one is hence far from being inappropriate.
\end{rmrk}

\subsection{Validity of the approximation for sufficiently decaying topographies}

\noindent Recalling that $(U_{app},N_{app})$ is given by
(\ref{ansatz}), we introduce the following quantities :
\begin{equation} \label{app}
v_{app}^{\;\;\varepsilon} = \frac{U_{app}+N_{app}}{2} \;\;,\;\;
\eta_{app}^{\;\;\varepsilon} = \frac{U_{app}-N_{app}}{2}\;,
\end{equation}
in order to prove the following proposition. 

\vspace{1em}

\begin{prop}\label{propo1}
Let $s \ge 2$, $\sigma \ge s+5$, $(v_0,\eta_0) \in H^{\sigma}(\R)$ and $b
\in W^{1,\infty}(\R)$. There exists $T_0 > 0$ and a unique family $(v^{\,\varepsilon}_{\Sigma},\eta^{\,\varepsilon}_{\Sigma})_{0\le\varepsilon\le\varepsilon_0} \in L^{\infty}([0,\frac{T_0}{\varepsilon}];H^{\sigma}(\R))$ of solutions of $(\Sigma)$ with initial data $(v_0,\eta_0)$. We define $(u_0,n_0) = (v_0+\eta_0,v_0-\eta_0)$. Then there exists a unique solution $(U_0,N_0)$ to the system $(\Sigma_{KdV})$ with initial data $(u_0,n_0)$ and this solution is bounded in $L^{\infty}([0,T_0];H^{\sigma}(\R))$.\\ Moreover, we have the following error estimate for all $t \in [0,\frac{T_0}{\varepsilon}]$ :
$$
\Big|(v^{\,\varepsilon}_{\Sigma},\eta^{\,\varepsilon}_{\Sigma})-(v^{\;\;\varepsilon}_{app},\eta^{\;\;\varepsilon}_{app})\Big|_{L^{\infty}([0,t];H^{s}(\R))} \le C(1+\left|(U_1,N_1)\right|_{L^{\infty}\left([0,T_0]\times[0,t];H^{s+3}(\R)\right)})\,\varepsilon^2 t\,,
$$
where $(v^{\;\;\varepsilon}_{app},\eta^{\;\;\varepsilon}_{app})$ is defined in (\ref{app}).
\end{prop}

\vspace{1em}

\begin{pf}
The result on the system $(\Sigma_{KdV})$ is a very classical result on the KdV equation that has been established f.e. by Bona and Smith (\cite{BS}, 1975) and we omit the proof here. \\
The leading terms $(U_0,N_0)$ and the correcting terms $(U_1,N_1)$ have been chosen such that $(U_{app},N_{app})$ is solution of the system $(\Gamma)$ with a residual of order $O(\varepsilon^2)$. This residual denoted by $(\varepsilon^2 R_1,\varepsilon^2 R_2)$ can be computed explicitly and we get :
$$
\begin{array}{l}
R_1 = -\frac{3}{4}\partial_x (U_0 U_1)-\frac{1}{4}\partial_x (N_0 N_1) - \frac{1}{4} \partial_x (U_0 N_1) - \frac{1}{4} \partial_x (U_1 N_0) + \frac{1}{2} b \partial_x U_1 \\ \hspace{2em}+ \frac{1}{4} (U_1+N_1) \partial_x b - a_1 \partial_x^3 U_1 + \frac{a_2+a_4}{2} \partial_x^2\partial_t U_1 + \frac{a_2-a_4}{2} \partial_x^2\partial_t N_1 
\\ \hspace{2em}
+ \frac{a_2+a_4}{2} \partial_x^2\partial_T U_0 + \frac{a_2-a_4}{2} \partial_x^2\partial_T N_0 - \varepsilon \left[ \frac{3}{8} \partial_x U_1^2 + \frac{1}{8} \partial_x N_1^2 + \frac{1}{4} \partial_x (U_1 N_1) \right]\,.
\end{array}
$$
with a similar expression holding for the residual $R_2$ of the second equation of $(\Gamma)$. We first use $(\Sigma_{KdV})$ and $(\Sigma_{corr})$ to express $\partial_T U_0$, $\partial_T N_0$, $\partial_t U_1$, $\partial_t N_1$ in terms of spatial derivatives of $U_0,N_0,U_1,N_1$ and then use the following standard estimates : one first applies the differentiel operator $\partial_x^s$ to each relation giving $R_1$ and $R_2$, then multiplies them by respectively $\partial_x^s R_1$ and $\partial_x^s R_2$, and finally integrates on $\R$. Using the fact that $s \ge 2$, this yields easily for all $t \in [0,\frac{T_0}{\varepsilon}]$ :
$$
\left|(R_1,R_2)\right|_{L^{\infty}\left([0,t];H^{s}(\R)\right)} \le C(1+\left|(U_1,N_1)\right|_{L^{\infty}\left([0,T_0]\times[0,t];H^{s+3}(\R)\right)})\;,
$$
with $C$ depending only on $|(U_0,N_0)|_{L^{\infty}\left([0,T_0];H^{s+5}(\R)\right)}$ and $|b|_{W^{1,\infty}(\R)}$. Inverting the diagonalization by plugging the relations $U_{app}=v_{app}+\eta_{app}$ and $N_{app} = v_{app}-\eta_{app}$ into $(\Gamma)$,
we easily deduce that $(v^{\;\;\varepsilon}_{app},\eta^{\;\;\varepsilon}_{app})$ is solution of the system $(\Sigma)$ with a residual bounded by $C\varepsilon^2(1+\left|(U_1,N_1)\right|_{L^{\infty}\left([0,T_0]\times[0,t];H^{s+3}(\R)\right)})$. Standard energy estimates - as described above - applied on the \emph{symmetric} Boussinesq system $(\Sigma)$ yield :
$$
|(v^{\,\varepsilon}_{\Sigma},\eta^{\,\varepsilon}_{\Sigma})-(v^{\;\;\varepsilon}_{app},\eta^{\;\;\varepsilon}_{app})|_{L^{\infty}([0,t];H^{s}(\R))} \le C(1+\left|(U_1,N_1)\right|_{L^{\infty}\left([0,T_0]\times[0,t];H^{s+3}(\R)\right)})\,\varepsilon^2 t\,,
$$ 
which ends the proof.
\end{pf}
An easy extension of this proposition is the following corollary which gives an error bound for the KdV approximation.  

\vspace{0.5em}

\begin{coro}\label{propo2}
Under the same hypothesis as in Proposition \ref{propo1}, we have the following error estimate for all $t \in [0,\frac{T_0}{\varepsilon}]$ :
$$
\Big|(v^{\,\varepsilon}_{\Sigma},\eta^{\,\varepsilon}_{\Sigma})-(v^{\;\;\;\varepsilon}_{KdV},\eta^{\;\;\;\varepsilon}_{KdV})\Big|_{L^{\infty}([0,t];H^{s}(\R))} \le C(1+\left|(U_1,N_1)\right|_{L^{\infty}\left([0,T_0]\times[0,t];H^{s+3}(\R)\right)})\,\varepsilon(1+\varepsilon t)\,,
$$
where $(v^{\;\;\;\varepsilon}_{KdV},\eta^{\;\;\;\varepsilon}_{KdV})$ is the uncoupled KdV approximation defined in (\ref{M}).
\end{coro}

\vspace{0.5em}

\begin{pf}
One has :
$$
\begin{array}{lll}
\vspace{0.5em}
(v^{\,\varepsilon}_{\Sigma},\eta^{\,\varepsilon}_{\Sigma})-(v^{\;\;\;\varepsilon}_{KdV},\eta^{\;\;\;\varepsilon}_{KdV})
& = &
(v^{\,\varepsilon}_{\Sigma},\eta^{\,\varepsilon}_{\Sigma}) - (v^{\;\;\varepsilon}_{app},\eta^{\;\;\varepsilon}_{app})
+ (v^{\;\;\varepsilon}_{app},\eta^{\;\;\varepsilon}_{app}) \\ & & - (v^{\;\;\;\varepsilon}_{KdV},\eta^{\;\;\;\varepsilon}_{KdV}) \\ \vspace{0.5em}
& = & (v^{\,\varepsilon}_{\Sigma},\eta^{\,\varepsilon}_{\Sigma})-(v^{\;\;\varepsilon}_{app},\eta^{\;\;\varepsilon}_{app}) + \varepsilon \Big(\frac{U_1+N_1}{2},\frac{U_1+N_1}{2}\Big)
\end{array}
$$
Using this relation and the error estimate coming from Proposition \ref{propo1} yields the result.
\end{pf}

\vspace{1em}

\noindent This corollary clearly states that the validity of the uncoupled KdV approximation only depends on the control of the correcting terms $(U_1,N_1)$ in $L^{\infty}\Big([0,T_0]\times[0,t]$; $H^{s+3}(\R)\Big)$ norm on the large time scale $1/\varepsilon$. From now on, these correctors become the center of our analysis.\\
As we saw earlier, the inhomogeneous transport equations that govern the evolution of the correctors
$(U_1,N_1)$ can be solved explicitly in terms of $U_0$ and $N_0$. Using the fact that the solution of the equation $(\partial_t+\partial_x) u = f$ is given by $u(t,x) = \int_0^t f(x-t+s) ds + u(0,x-t)$ and that $U_1(t=0)=0$, we
thus get the following expression :
\begin{eqnarray} \label{explicit}
U_1(T,t,x) & = & -\frac{1}{16}(N_0^2(T,x+t)-N_0^2(T,x-t)) + \frac{a_2-a_4}{4}(\partial_x^2 N_0(T,x+t) \nonumber \\ 
& & - \partial_x^2 N_0(T,x-t)) - \frac{1}{8} U_0(T,x-t)(N_0(T,x+t)-N_0(T,x-t)) \nonumber \\
& & + \frac{1}{4} U_0(T,x-t) (b(x)-b(x-t)) 
- \frac{1}{4}\partial_x U_0(T,x-t) \nonumber \\ 
& & \int_0^t{N_0(T,x-t+2s)ds} + \frac{1}{2}\partial_x U_0(T,x-t)\int_0^t{b(x-t+s)ds} 
\nonumber \\ 
& & + \frac{1}{4}\int_0^t{\partial_x b(x-t+s)N_0(T,x-t+2s)ds}\;, 
\end{eqnarray}
and a similar expression holds for $N_1$ :
\begin{eqnarray} \label{explicit2}
N_1(T,t,x) & = & -\frac{1}{16}(U_0^2(T,x-t)-U_0^2(T,x+t)) - \frac{a_2-a_4}{4}(\partial_x^2 U_0(T,x-t) \nonumber \\ 
& & - \partial_x^2 U_0(T,x+t)) - \frac{1}{8} N_0(T,x+t) (U_0(T,x-t)-U_0(T,x+t)) \nonumber \\
& & - \frac{1}{4} N_0(T,x+t) (b(x)-b(x+t)) 
- \frac{1}{4}\partial_x N_0(T,x+t) \nonumber \\ 
& & \int_0^t{U_0(T,x+t-2s)ds} - \frac{1}{2}\partial_x N_0(T,x+t)\int_0^t{b(x+t-s)ds} 
\nonumber \\ 
& & - \frac{1}{4}\int_0^t{\partial_x b(x+t-s)U_0(T,x+t-2s)ds}\;, 
\end{eqnarray}
\noindent We here only deal with the
case of $U_1$ since all the method can easily be adapted to the case
of $N_1$.

\vspace{1em}

\noindent The corrector $U_1$ is analysed in the following way : let $T_0 \ge 0$,
$s \ge 2$, $\sigma \ge s+5$ and $(u_0,n_0)$ in
$H^{\sigma}(\R)$. We know that the solutions $(U_0,N_0)$ of the problem
$(\Sigma_{KdV})$ with initial data $(u_0,n_0)$ are bounded in 
$L^{\infty}([0,T_0];H^{\sigma}(\R))^2$. We suppose here that the
bottom topography $b$ is bounded in $W^{1,\infty}(\R)$. 
Under these circumstances, it clearly appears that the first
four terms of the expression of $U_1$ are bounded in
$L^{\infty}([0,T_0]\times[0,t];H^{s+3}(\R))$. Only the last four terms can
be problematic and deserve a precise treatment. \\
If $(U_0,N_0)$ and $b$ come with a sufficient decay
rate at infinity, we can straightforwardly
control these terms. To this end and following \cite{SW}, we introduce the following weighted
Sobolev space $H^{s,\alpha}$ for all $s \in \N$ and $\alpha \in \R$ :
$$
H^{s,\alpha} = \Big\{u \in H^s(\R) \,/\, |u|_{H^{s,\alpha}}^2 \equiv
  \sum_{k=0}^{s} \Big|(1+x^2)^{\alpha/2}\,\frac{\partial^k
      u}{\partial x^k}\Big|_{L^2(\R)} < \infty \Big\}\;.
$$
We can now state our first theorem on the validity of the
approximation for sufficiently decaying initial data and bottom topography.

\vspace{1em}

\begin{thrm} \label{thrm1}
Let $s \ge 2$, $\sigma \ge s+5$, $(v_0,\eta_0) \in
H^{\sigma,1}(\R)^2$ and $b \in H^{s+4,1}(\R)$. There exists $T_0 > 0$ and a unique family $(v^{\,\varepsilon}_{\Sigma},\eta^{\,\varepsilon}_{\Sigma})_{0\le\varepsilon\le\varepsilon_0} \in L^{\infty}([0,\frac{T_0}{\varepsilon}];H^{\sigma}(\R))$ of solutions of $(\Sigma)$ with initial data $(v_0,\eta_0)$. We define $(u_0,n_0) = (v_0+\eta_0,v_0-\eta_0)$. Then the solution $(U_0,N_0)$ of the system $(\Sigma_{KdV})$ with initial data $(u_0,n_0)$ is bounded in $L^{\infty}([0,T_0];H^{\sigma,1}(\R))$. Moreover, we have the following error estimate for all $t \in [0,\frac{T_0}{\varepsilon}]$ :
$$
\Big|(v^{\,\varepsilon}_{\Sigma},\eta^{\,\varepsilon}_{\Sigma})-(v^{\;\;\;\varepsilon}_{KdV},\eta^{\;\;\;\varepsilon}_{KdV})\Big|_{L^{\infty}([0,t];H^{s}(\R))} \le C\varepsilon(1+\varepsilon t)\,,
$$
where $(v^{\;\;\;\varepsilon}_{KdV},\eta^{\;\;\;\varepsilon}_{KdV})$ is the uncoupled KdV approximation defined in (\ref{M}).
\end{thrm}

\vspace{1em}

\begin{pf}
We know from \cite{Kato} and \cite{SW} that the KdV equation propagates the regularity of initial data taken in weighted Sobolev spaces and we omit the proof here. The end of the proof is devoted to the estimate of $|(U_1,N_1)|_{L^{\infty}([0,T_0]\times [0,t];H^{s+3}(\R])^2}$.
The work of Lannes in \cite{Lannes3} is here very useful to control this quantity. Indeed, using the equations $(\Sigma_{corr})$, the fact that $U_0(T,.)$, $N_0(T,.)$ and $b$ are bounded in $H^{\sigma,1}(\R)$, and Proposition 3.5 of \cite{Lannes3}, one finally obtains the estimate :
$$
|(U_1,N_1)|_{L^{\infty}([0,T_0]\times [0,t];H^{s+3}(\R])^2} \le C(|b|_{H^{\sigma,1}(\R)},(U_0,N_0)_{L^{\infty}([0,t];H^{\sigma,1}(\R))^2})\;.
$$
Pluging this last estimate into the result of Proposition 1.1 ends the proof.
\end{pf}

\vspace{1em}

\begin{rmrk}
As specified in \cite{BCL}, this approximation diverges on a large
time scale in the periodic framework unless we specify a zero mass
assumption on the initial data $u_0$ and $n_0$. This drawback is dealt
with at the end of the next section. Until then, the results provided can be extended to the periodic framework with this zero mass hypothesis.
\end{rmrk}
\begin{rmrk}
It is worth pointing out that the validity of the uncoupled KdV
approximation for the Boussinesq system $(\Sigma)$ is enough to
demonstrate its validity regarding the water waves problem. Indeed, we
can deduce from \cite{FC} that the error estimate between the solutions $(v^{\,\varepsilon}_{\Sigma},\eta^{\,\varepsilon}_{\Sigma})$ of $(\Sigma)$ and the solutions of the water waves problem is of order $O(\varepsilon(1+\varepsilon t))$. An error estimate between the solutions of the water waves problem and the KdV approximation $(v^{\;\;\;\varepsilon}_{KdV},\eta^{\;\;\;\varepsilon}_{KdV})$ can thus be immediately deduced from the results of this paper.
\end{rmrk}

\section{A topographically modified KdV approximation}

In this section we discuss the validity of the previously derived
uncoupled KdV approximation on a large time scale for different
bottom topographies. We demonstrate its validity for less restrictive bottoms,
but provide two examples of simple bottoms for which the approximation
diverges. A new approximation that takes the bottom into account
is finally derived.

\subsection{Discussion on the validity of the approximation}

\noindent Starting from the previous theorem, it is worth wondering if
this one holds for less restrictive initial data and bottoms,
i.e. without any condition of a sufficient decay rate at infinity. In
this view,  we focus in a more general way on the last three terms of
$U_1$ by supposing that $(u_0,n_0)$ is bounded in
$L^{\infty}([0,t];H^{\sigma}(\R))^2$, which is propagated by the KdV
equation on $(U_0,N_0)$ (see \cite{KPV}).
Using Cauchy-Schwarz inequality on the first two terms
and Proposition 3.2 of \cite{Lannes3} on the last term, we can write the
following controls for all $t \in [0,\frac{T_0}{\varepsilon}]$, $s \ge
2$ and $\sigma \ge s+5$:

$$
\begin{array}{l}
\vspace{1em}
\Big|\partial_x U_0(T,.-t) \int_0^t N_0(T,.-t+2s) ds\Big|_{H^s(\R)} \le C_1 \sqrt{t}\;,\\
\vspace{1em}
\Big|\partial_x U_0(T,.-t) \int_0^t b(.-t+s) ds\Big|_{H^s(\R)} \le C_2 \left|b\right|_{L^2(\R)} \sqrt{t} \;,\\
\Big|\int_0^t \partial_x b(.-t+s) N_0(T,.-t+2s) ds\Big|_{H^s(\R)} \le C_3 \left|\partial_x b\right|_{H^s(\R)} \sqrt{t}\;,
\end{array}
$$
where the constants $C_1,C_2,C_3$ depend exclusively on $|(U_0,N_0)|_{L^{\infty}([0,t];H^{\sigma}(\R))^2}$.\\
These preliminary estimates are at the heart of the proof of the
following theorem.

\vspace{1em}

\begin{thrm} \label{thrm2}
Let $s \ge 2$, $\sigma \ge s+5$, $(v_0,\eta_0) \in
H^{\sigma}(\R)^2$, $b \in H^{s+4}(\R)$. There exists $T_0 > 0$ and a unique family $(v^{\,\varepsilon}_{\Sigma},\eta^{\,\varepsilon}_{\Sigma})_{0\le\varepsilon\le\varepsilon_0} \in L^{\infty}([0,\frac{T_0}{\varepsilon}];H^{\sigma}(\R))$ of solutions of $(\Sigma)$ with initial data $(v_0,\eta_0)$. We define $(u_0,n_0) = (v_0+\eta_0,v_0-\eta_0)$. Then
the solution $(U_0,N_0)$ of the system $(\Sigma_{KdV})$ with initial
data $(u_0,n_0)$ is bounded in
$L^{\infty}([0,T_0];H^{\sigma}(\R))$. Moreover, we have the following error estimate for all $t \in [0,\frac{T_0}{\varepsilon}]$ :
$$
\Big|(v^{\,\varepsilon}_{\Sigma},\eta^{\,\varepsilon}_{\Sigma})-(v^{\;\;\;\varepsilon}_{KdV},\eta^{\;\;\;\varepsilon}_{KdV})\Big|_{L^{\infty}([0,t];H^{s}(\R))} \le C\varepsilon \sqrt{t} (1+\varepsilon t)\,,
$$
where $(v^{\;\;\;\varepsilon}_{KdV},\eta^{\;\;\;\varepsilon}_{KdV})$ are as defined in (\ref{M}).
\end{thrm}

\vspace{1em}

\begin{pf}
Using the three previous inequalities, one obtains :
$$
|(U_1,N_1)|_{L^{\infty}([0,T_0]\times [0,t];H^{s+3}(\R))} \le C \sqrt{t}\;.
$$
where $C=C(|b|_{H^{s+4}(\R)},(U_0,N_0)_{L^{\infty}([0,t];H^{\sigma}(\R))^2})$.
The final result follows from Corollary \ref{propo2}.
\end{pf}

\vspace{1em}

\noindent This theorem proves that the approximation is less
precise on a large time scale if we remove the
assumption of a sufficient decay rate at infinity. And yet, it is worth 
pointing out that the regularity imposed on $b$ in this theorem
excludes many physical cases of interest. We focus from now on two simple examples of bottoms which do not fall into the scope of Theorem \ref{thrm2} : a regular step,
and a slowly varying sinusoidal bottom. Our goal is to emphasize the
fact that the approximation
$(v^{\;\;\;\varepsilon}_{KdV},\eta^{\;\;\;\varepsilon}_{KdV})$ diverges from the
exact solution
$(v^{\,\varepsilon}_{\Sigma},\eta^{\,\varepsilon}_{\Sigma})$ in these two
simple cases. To deal with such bathymetries, a topographically modified KdV approximation is derived at the end of the section.

\vspace{1em}

\noindent In order to simplify the analysis, we only consider the approximation corresponding to $a_1 = 1/6,a_2=0,a_4=0$ which is obtained for $\theta=\sqrt{2/3},\lambda_1=1,\lambda_2=1$, and the case of a wave propagating to the right. This last condition is realized by taking $n_0=0$, which implies that $N_0=N=0$.

\subsubsection{The case of a step}

\noindent We consider here a bottom whose shape corresponds to a regular step. The interest of such an example is that in this case, $b \notin L^2(\R)$. \\
The bottom is defined as follows :
\begin{equation}\label{step}
b(x) = \;
\left\{
\begin{array}{l}
\vspace{0.5em}
0\;,\;\;\forall x \le 0 \;,\\
\vspace{0.5em}
\frac{A}{2}\left(1+\sin\Big(\frac{\pi}{l}(x-\frac{l}{2})\Big)\right) \;,\;\;\forall x \in [0,l] \;,\\
A\;,\;\;\forall x \ge l \;.
\end{array}
\right.
\end{equation}
For a right going wave, the system $(\Sigma_{KdV})$ is reduced to the simple KdV equation :
$$
\partial_T U_0 + \frac{3}{8} \partial_x U_0^2 + \frac{1}{6}
\partial_x^3 U_0 = 0\;,
$$
and we chose the initial condition $u_0$ such that the solution of this equation is a positive soliton which propagates to the right. \\
We write the explicit expression of the corrector $U_1$ when $N_0=0$ :
\begin{equation*}
\begin{array}{l}
U_1(t,x) =  \frac{1}{4}U_0(x-t)(b(x)-b(x-t)) + \frac{1}{2}\partial_x U_0(x-t)\int_0^t{b(x-t+s)ds}\;.
\end{array}
\end{equation*}
In this expression, 
the only possibly secularly growing term is $\partial_x U_0(T,x-t) \int_0^t
b(x-t+s) ds$. The time evolution in amplitude of this
term is obviously led by the evolution of $\int_0^t b(x-t+s)$ for all
$x \in \R$. When the bottom is a step as defined in (\ref{step}), this
integral essentially grows linearly in time. We now prove that because
of this, $|U_1|_{L^{\infty}([0,T_0]\times [0,t];H^{s+3}(\R))}$ grows
linearly in time. Let $s \ge 2$ and $\sigma \ge s+5$. Starting from the expression of $U_1$, we get for all $t \in [0,\frac{T_0}{\varepsilon}]$  the following estimates :
\begin{eqnarray*}
|U_1(T,t,\cdot)|_{H^{s+3}(\R)}
& \ge & \Big|\frac{1}{2} \partial_x U_0(T,\cdot-t) \int_0^t b(\cdot-t+s)
ds\Big|_{H^{s+3}(\R)} - C \;,\nonumber \\
\end{eqnarray*}
with $C = \Big|\frac{1}{4}U_0(T,\cdot-t)(b(\cdot)-b(\cdot-t))\Big|_{H^{s+3}}
\le \frac{1}{2} |b|_{L^{\infty}} |U_0|_{L^{\infty}([0,t];H^{s+3})}
\equiv C_0$,
\begin{eqnarray*}
|U_1(T,t,\cdot)|_{H^{s+3}(\R)}
& \ge & \frac{1}{2} |\partial_x U_0(T,\cdot-t) \int_0^t b(\cdot-t+s)
ds|_{L^2(\R)} - C_0 \;,\\
& = & \frac{1}{2} \sqrt{\int_0^{\infty} |\partial_x U_0(T,x-t)|^2
    \; \Big|\int_0^t b(x-t+s)ds\Big|^2 dx} - C_0 \;,
\end{eqnarray*}
since $\int_0^t b(x-t+s)ds = 0 \;,\;\forall x \le 0 \;$,
\begin{eqnarray*}
\hspace{8em} & \ge & \frac{1}{2} \sqrt{\int_{l+t}^{\infty} |\partial_x
  U_0(T,x-t)|^2 \; \Big|\int_{x-t}^{x} b(s)ds\Big|^2 dx} - C_0 \;, \\
& = & \frac{1}{2} A t \, \sqrt{\int_{l+t}^{\infty} |\partial_x U_0(T,x-t)|^2 dx} - C_0 \;,
\end{eqnarray*}
since $\int_{x-t}^x b(s)ds = A t \;,\;\forall x \ge l+t \;$,
\begin{eqnarray*}
\hspace{8em} & = & \frac{1}{2} A\,t\, \sqrt{\int_{l}^{\infty} |\partial_x U_0(T,x)|^2 dx} - C_0 \;,
\end{eqnarray*}
which implies that
\begin{equation}  \label{1}
\left|U_1\right|_{L^{\infty}([0,T_0]\times [0,t];H^{s+3}(\R))} \ge C_1 t - C_0\;,
\end{equation}
where the last positive constant $C_1$ only depends on $|\partial_x U_0|_{L^2(\R)}$.

\vspace{1em}

\noindent This linear growth of $|U_1|_{L^{\infty}([0,T_0]\times
  [0,t];H^{s+3}(\R))}$ is sharp since it follows from the explicit
expression of $U_1$ that this growth is at most linear.
Furthermore, we recall that
\begin{equation} \label{3}
(v^{\,\varepsilon}_{\Sigma},\eta^{\,\varepsilon}_{\Sigma})-(v^{\;\;\;\varepsilon}_{KdV},\eta^{\;\;\;\varepsilon}_{KdV}) =  (v^{\,\varepsilon}_{\Sigma},\eta^{\,\varepsilon}_{\Sigma})-(v^{\;\;\varepsilon}_{app},\eta^{\;\;\varepsilon}_{app}) + \varepsilon \left(\frac{U_1+N_1}{2},\frac{U_1-N_1}{2}\right)\;.
\end{equation}
Using this relation, (\ref{1}) and Proposition \ref{propo1}, we get that there exists two constants $C_2$ and $C_3$ and a time $T_1$ independent of $\varepsilon$ such that $\forall t \in [T_1,\frac{T_0}{\varepsilon}]$,
$$
\Big|(v^{\,\varepsilon}_{\Sigma},\eta^{\,\varepsilon}_{\Sigma})-(v^{\;\;\;\varepsilon}_{KdV},\eta^{\;\;\;\varepsilon}_{KdV})\Big|_{L^{\infty}([0,t];H^{s}(\R))}
\ge \left|C_2 (1+t) \varepsilon - C_3 (1+t)\varepsilon^2 t\right|  \;.
$$
We finally deduce that there exists two constants $C$ and $C'$ such
that $\forall t \in [T_1,\frac{T_0}{\varepsilon}]$,
$$
\Big|(v^{\,\varepsilon}_{\Sigma},\eta^{\,\varepsilon}_{\Sigma})-(v^{\;\;\;\varepsilon}_{KdV},\eta^{\;\;\;\varepsilon}_{KdV})\Big|_{L^{\infty}([0,t];H^{s}(\R))}
\ge  C \varepsilon t |C'-\varepsilon t|\;.
$$
This proves that in this case, the error is of order $O(1)$ on
times of order $O(1/\varepsilon)$, and the usual KdV approximation is
not valid for such a topography.

\subsubsection{The case of a sinusoidal bottom}

\noindent We consider here a bottom defined as follows :
\begin{equation}\label{sinus}
b(x) = A \sin(\varepsilon x)\;,\;\;\forall x \in \R \;.
\end{equation}
where $U_0$ is again a soliton propagating to the right.\\
We mention that such a type of periodic bottom varying on a slow
spatial scale has been studied in \cite{Sulem} by
Craig-Guyenne-Nicholls-Sulem, with the difference that the authors
authorized the bottom to vary also on a small spatial scale. \\
Again, the amplitude of the term $\partial_x U_0(T,x-t) \int_0^t
b(x-t+s) ds$ evolves in time according to $\int_0^t b(x-t+s) ds$. Let
us have a look at this quantity for all $x \in \R$ and $t \ge 0$ :
\begin{eqnarray*}
\int_0^t b(x-t+s) ds & = & \int_{x-t}^x \sin(\varepsilon x) ds  \\
& = & -\frac{A}{\varepsilon} \left[\cos(\varepsilon x) -
  \cos(\varepsilon (x-t))  \right] \\
& = & \frac{2A}{\varepsilon} \sin\left(\varepsilon (x-\frac{t}{2})\right)
\sin\left(\frac{\varepsilon t}{2}\right)\;.
\end{eqnarray*}
We can see that the amplitude of this term is of order
$O(1/\varepsilon)$. We now demonstrate that it is also the case for
the corrector $U_1$ :
\begin{eqnarray*} \label{f}
\left|U_1(T,t,\cdot)\right|_{H^{s+3}(\R)}
& \ge & \left|\frac{1}{2} \partial_x U_0(T,\cdot-t) \int_0^t b(\cdot-t+s) ds\right|_{H^{s+3}(\R)}
- C_0\;,\nonumber \\
\end{eqnarray*}
\begin{eqnarray}
& \ge & \frac{1}{2} \left|\partial_x U_0(T,\cdot-t) \int_0^t b(\cdot-t+s)
  ds\right|_{L^2(\R)} -C_0 \;,\nonumber \\
& = & \frac{A}{\varepsilon} \sqrt{\int_{-\infty}^{\infty} \left|\partial_x
    U_0(T,x-t)\right|^2 \sin^2\left(\varepsilon (x-\frac{t}{2})\right)
\sin^2\left(\frac{\varepsilon t}{2}\right) dx} - C_0
\;,\nonumber \\
& = & \frac{A}{\varepsilon} |\sin^2\left(\frac{\varepsilon t}{2}\right)| \sqrt{\int_{-\infty}^{\infty} \left|\partial_x
    U_0(T,x-t)\right|^2 \sin^2\left(\varepsilon (x-\frac{t}{2})\right) dx} - C_0
\;.\nonumber \\
& & 
\end{eqnarray}
At this point, we remark that
$$
0 \le \int_{-\infty}^{\infty} \left|\partial_x
    U_0(T,x-t)\right|^2 \sin^2\left(\varepsilon (x-\frac{t}{2})\right)
  dx \le \int_{-\infty}^{\infty} \left|\partial_x
    U_0(T,x-t)\right|^2 dx\;,
$$
We hence deduce that for all $t \ge 0$ there exists $\alpha(t) \in \R$ such that
$$
 \int_{-\infty}^{\infty} \left|\partial_x
    U_0(T,x-t)\right|^2 \sin^2\left(\varepsilon (x-\frac{t}{2})\right)
  dx = \sin^2\left(\alpha(t)\right) \int_{-\infty}^{\infty} \left|\partial_x
    U_0(T,x-t)\right|^2 dx\;.
$$
Pluging this one into (\ref{f}) leads to
$$
\left|U_1(T,t,\cdot)\right|_{H^{s+3}(\R)} \ge \frac{A}{\varepsilon}
\left|\sin^2\left(\frac{\varepsilon
      t}{2}\right)\,\sin^2\left(\alpha(t)\right)\right| \Big|\partial_x
U_0(T,\cdot)\Big|_{L^2(\R)} - C_0\;,
$$
which finally implies that there exists a constant $C_1$ such that
$$
|U_1|_{L^{\infty}([0,T_0]\times [0,t];H^{s+3}(\R))} \ge \frac{C_1}{\varepsilon}-C_0\;.
$$
\noindent Using this result and the same technique as in the previous example leads to the same conclusion : 
the uncoupled KdV approximation diverges on a large time scale in this case too.

\subsection{A topographically modified approximation}

Both examples clearly show the
invalidity of the approximation on a large time scale if we consider
general bottoms topographies $b$ which do not have specific decay
properties at infinity. Therefore, we need to modify the
usual KdV approximation to be able to handle general bathymetries. \\
All the previous analysis has shown that two terms of the r.h.s. of the explicit expression (\ref{explicit}) of $U_1$ may
exhibit a secular growth, and cause the approximation to diverge on a long time scale: 
these are $\frac{1}{2} \partial_x U_0(T,x-t)
\int_0^t b(x-t+s) ds$ and $\frac{1}{4} \int_0^t \partial_x b(x-t+s)
N_0(T,x-t+2s) ds$. As far as the expression (\ref{explicit2}) of $N_1$ is concerned, the same possibly
problematic terms are $\frac{-1}{2} \partial_x N_0(T,x+t) $\\$ \int_0^t
b(x+t-s) ds$ and $\frac{-1}{4} \int_0^t \partial_x b(x+t-s)
U_0(T,x+t-2s) ds$. The idea is as follows: rather than treating these
terms as correcting terms - which invalidates the approximation for general bathymetries - we can include them with the leading order one terms $U_0$ and $N_0$ in the final approximation.

\vspace{1em}

\noindent This idea leads us to propose the following topographically
modified KdV approximation which is an alternative version of
$(\mathcal{M})$ :

\begin{eqnarray} \label{MB}
(\mathcal{M}_b)\,\left\{
\begin{array}{rrl} 
v_{KdV}^{\;\,\varepsilon,b} & = & \frac{U_0+N_0}{2} +
\frac{\varepsilon}{4} \left[\partial_x U_0(T,x-t) \int_0^t b(x-t+s) ds
\right. \\ & & \left.
\hspace{6em} - \partial_x N_0(T,x+t) 
\int_0^t b(x+t-s) ds \right. \\ & & \left. 
\hspace{6em} + \frac{1}{2} \int_0^t \partial_x b(x-t+s) N_0(T,x-t+2s) ds 
\right. \\ & & \left. 
\hspace{6em} - \frac{1}{2} \int_0^t \partial_x b(x+t-s) U_0(T,x+t-2s) ds 
\right. \\ & & \left. 
\hspace{6em} + \frac{1}{2} U_0(T,x-t)\left(b(x)-b(x-t)\right)
\right. \\ & & \vspace{1em} \left. 
\hspace{6em} + \frac{1}{2} N_0(T,x+t)\left(b(x+t)-b(x)\right) \right] \;\;, \\ 
\eta_{KdV}^{\;\,\varepsilon,b} & = & \frac{U_0-N_0}{2} +
\frac{\varepsilon}{4} \left[\partial_x U_0(T,x-t) \int_0^t b(x-t+s) ds
\right. \\ & & \left.
\hspace{6em} - \partial_x N_0(T,x+t) 
\int_0^t b(x+t-s) ds \right. \\ & & \left. 
\hspace{6em} + \frac{1}{2} \int_0^t \partial_x b(x-t+s) N_0(T,x-t+2s) ds 
\right. \\ & & \left. 
\hspace{6em} - \frac{1}{2} \int_0^t \partial_x b(x+t-s) U_0(T,x+t-2s) ds 
\right. \\ & & \left. 
\hspace{6em} + \frac{1}{2} U_0(T,x-t)\left(b(x)-b(x-t)\right)
\right. \\ & & \vspace{1em} \left. 
\hspace{6em} + \frac{1}{2} N_0(T,x+t)\left(b(x+t)-b(x)\right) \right] \;\;.
\end{array}
\right.
\end{eqnarray}
where $U_0$ and $N_0$ are still solutions of the system $(\Sigma_{KdV}^{\;\;\;\varepsilon})$, and where all the possibly secularly growing terms of (\ref{explicit}) and (\ref{explicit2}) have been included with the order one terms. The physical role of each of these additionnal topography-dependent terms is discussed in the last section, where we validate this model numerically on the previous examples of a step and a sinusoidal bathymetry.

\vspace{1em}

\begin{rmrk}
We have here also included the terms
$U_0(T,x-t)(b(x)-b(x-t))$ and
$N_0(T,x+t)\left(b(x+t)-b(x)\right)$ even if these terms remain
bounded 
indepently of $\varepsilon$ for all time. The reason of this choice is
that we are interested in their physical meaning. Indeed, we further
see - in the last section - that they are responsible
for the reproduction of the phenomenon of shoaling. We hence decided
to include these terms in the approximation. 
\end{rmrk}

\vspace{1em}

\noindent The main advantage of this modification relies in the following remark : now that the bottom terms have been included with the leading order terms in the approximation, we can easily see that the correcting terms $U_1$ and $N_1$ solve a different equation. Indeed, the equations on $U_1$ and $N_1$ become :
\begin{equation}
(\Sigma_{corr}^{\;\;b})\;\left\{
\begin{array}{l}
\vspace{0.5em}
(\partial_t + \partial_x) U_1 = - \frac{1}{8} \partial_x N_0^2 - \frac{1}{4} \partial_x (U_0 N_0) + \frac{a_2-a_4}{2}\, \partial_x^3 N_0\;, \\
(\partial_t - \partial_x) N_1 = - \frac{1}{8} \partial_x U_0^2 - \frac{1}{4} \partial_x (U_0 N_0) - \frac{a_2-a_4}{2}\, \partial_x^3 U_0\;.
\end{array}
\right.
\end{equation}
It is clear here that all the possibly secularly growing terms of the
correctors $(U_1,N_1)$ have been removed.

\vspace{0.5em}

\begin{rmrk}
Numerically speaking, this modified version is quite interesting since the topographical terms are computed explicitly from the solution
of the KdV equations. We thus expect the numerical simulation of this
model to be faster than the one of the symmetric Boussinesq model
$(\Sigma)$. This point is checked in the last section. 
\end{rmrk}

\vspace{0.5em}

\noindent In the periodic framework, we saw that the usual approximation is not valid on a large time scale because of the linear growth in time of the term $\partial_x U_0(T,x-t)\int_0^t N_0(T,x-t+2s)ds$ in $U_1$, unless we specify a zero mass assumption on the initial data $u_0$ and $n_0$. Once more, we can propose a valid approximation just by including this term in the order one terms of the ansatz. We conclude this section with the proposition of a new approximation in the periodic framework :
\begin{eqnarray} \label{newp}
(\mathcal{M}_{\;\;b}^{per})\;
\left\{
\begin{array}{lll}
v_{\,KdV}^{\varepsilon,b\,per} & = &  v_{KdV}^{\;\,\varepsilon,b} - \frac{\varepsilon}{8} \left[ \partial_x U_0(T,x-t)\int_0^t N_0(T,x-t+2s)ds 
\right. \\ & & \vspace{1em} \left.
+ \partial_x N_0(T,x+t)\int_0^t U_0(T,x+t-2s)ds\right] \;\;, \\ 
\eta_{\,KdV}^{\varepsilon,b\,per} & = & \eta_{KdV}^{\;\,\varepsilon,b} - \frac{\varepsilon}{8} \left[ \partial_x U_0(T,x-t)\int_0^t N_0(T,x-t+2s)ds 
\right. \\ & & \vspace{1em} \left.
- \partial_x N_0(T,x+t)\int_0^t U_0(T,x+t-2s)ds\right]\;\;.
\end{array}
\right.
\end{eqnarray}

\vspace{1em}

To end this section, we would like to add a few words on the works - referenced in the introduction - \cite{Ostrovskii,Kakutani,Johnson,Miles,Newell} and \cite{Iguchi}. These papers are mainly devoted to the derivation of KdV equations in a variable medium, but differ from the present work in the way the topography is accounted for. Indeed, the bathymetry does not show up as additionnal terms computed from the solutions of the uncoupled KdV equations like here, but as variable topography-dependent coefficients in the KdV equations themselves. A good example of this result is the coupled KdV model obtained by Iguchi in \cite{Iguchi}. Their approach is thus very different from ours, and lead to additionnal difficulties such as the well-posedness of these topography-dependent KdV-like equations. However, Iguchi \cite{Iguchi} has been able to show that his coupled KdV approximation is well-posed in some Sobolev spaces and valid - in the meaning of the precision of the approximation given by the solutions - on a long time scale, as long as - among other conditions - the topography $b$ belongs to $W^{k,\infty}(\R)$ for $k$ large enough. This model is hence quite pertinent in the present context of realistic bathymetries and offers a good alternative to our topographically modified KdV approximation.

\noindent Nevertheless, it is worth pointing out that there actually exists a way to recover variable-coefficient KdV equations as in \cite{Ostrovskii,Kakutani,Johnson,Miles,Newell,Iguchi}, simply by modifying the ansatz (\ref{ansatz}) as follows :
\begin{eqnarray} \label{ansatz2}
\left\{
\begin{array}{l}
\vspace{1em}
U_{app}(t,x) = (1+\frac{\varepsilon}{4}\,b)\,U_0\left(T, x-(1-\frac{\varepsilon}{2}\,b) t\right)+\varepsilon U_1(T,t,x)\;, \\
N_{app}(t,x) = (1+\frac{\varepsilon}{4}\,b)\,N_0\left(T, x+(1-\frac{\varepsilon}{2}\,b) t\right)+\varepsilon N_1(T,t,x)\;,
\end{array}
\right.
\end{eqnarray}
Plugging this ansatz into $(\Gamma)$, neglecting all the $O(\varepsilon^2)$ terms but the $\varepsilon^2 b \partial_x U_1$ and $\varepsilon^2 b \partial_x N_1$ ones, and then proceeding as in section 1 leads finally to the same KdV equations on $U_0$ and $N_0$, but to a new system on the correcting terms :
\begin{equation}
(\Sigma_{corr}^{\;\;2})\;\left\{
\begin{array}{l}
\vspace{0.5em}
(\partial_t + (1-\frac{\varepsilon}{2}\,b)\partial_x) U_1 = - \frac{1}{8} \partial_x N_0^2 - \frac{1}{4} \partial_x (U_0 N_0) + \frac{1}{4} N_0 \partial_x b + \frac{a_2-a_4}{2}\, \partial_x^3 N_0\;, \\
(\partial_t - (1-\frac{\varepsilon}{2}\,b)\partial_x) N_1 = - \frac{1}{8} \partial_x U_0^2 - \frac{1}{4} \partial_x (U_0 N_0)  - \frac{1}{4} U_0 \partial_x b - \frac{a_2-a_4}{2}\, \partial_x^3 U_0\;.
\end{array}
\right.
\end{equation}
Finally, defining :
$$
\left\{
\begin{array}{l}
\vspace{0.25em}
u(t,x) = (1+\frac{\varepsilon}{4}\,b)\,U_0\left(T,x-(1-\frac{\varepsilon}{2}\,b) t\right)\;,\\
n(t,x) = (1+\frac{\varepsilon}{4}\,b)\,N_0\left(T,x+(1-\frac{\varepsilon}{2}\,b) t\right)\;,
\end{array}
\right.
$$
yields the final set of uncoupled KdV equations :
\begin{equation}
(\Sigma_{KdV}^{\;\;2})\;\left\{
\begin{array}{l}
\vspace{0.5em}
\partial_t u + \partial_x u + \varepsilon \Big[ \frac{3}{4} u \partial_x u + \frac{1}{6} \partial_x^3 u - \frac{1}{2}b \partial_x u - \frac{1}{4}\partial_x b\,u  \Big] = 0\;,\\
\partial_t n - \partial_x n + \varepsilon \Big[ \frac{3}{4} n \partial_x n - \frac{1}{6} \partial_x^3 n + \frac{1}{2}b \partial_x n + \frac{1}{4}\partial_x b\,n  \Big] = 0\;,
\end{array}
\right.
\end{equation}
and the final KdV approximation 
\begin{equation} \label{MB2}
(\mathcal{M}_b^2)\,\left\{
\begin{array}{rrl} 
v_{KdV}^{\;\,\varepsilon,b} & = & \frac{u+n}{2} +
\frac{\varepsilon}{8} \int_0^t \partial_x b(x-P_b^{\varepsilon}(t-s)) N_0(T,x+s-P_b^{\varepsilon}(t-s)) ds  \;\;, \\ 
\eta_{KdV}^{\;\,\varepsilon,b} & = & \frac{u-n}{2} -
\frac{\varepsilon}{8} \int_0^t \partial_x b(x+P_b^{\varepsilon}(t-s)) U_0(T,x-s+P_b^{\varepsilon}(t-s)) ds  \;\;, \\
\end{array}
\right.
\end{equation}
with $P_b^{\varepsilon} = 1-\frac{\varepsilon}{2}\,b$, and where we need to include the possibly secularly growing terms in the approximation. This time, the derived KdV equations are topography-dependent and look very similar to the one derived in f.e. \cite{Ostrovskii,Kakutani,Johnson}. This alternative formulation remains uncoupled and allow us to establish a link between our approach and the older works on the KdV approximation over uneven bottoms.

\section{Numerical comparison of the models}

This section is devoted to the numerical comparison of the different
models involved in this article. We compare here three models : the
symmetric Boussinesq system $(\Sigma)$ coming from \cite{FC}, the
usual uncoupled KdV approximation justified by Schneider-Wayne (\cite{SW}, flat bottoms) and Iguchi (\cite{Iguchi}, uneven
bottoms), and finally the topographically modified KdV
approximation. The aim is here to compare these three models for two
non trivial examples of topography : a step and a slowly varying
sinusoidal bottom.

\subsection{Numerical schemes}

Our goal is to compare three models, the symmetric Boussinesq one, the
usual KdV approximation $(\mathcal{M})$ and its topographically modified version
$(\mathcal{M}_b)$. The comparison is made for a solitary wave propagating to the right above two
topographies :  a step and a slowly varying
sinusoidal bottom. 
We use for the Boussinesq system $(\Sigma)$ and the KdV equations $(\Sigma_{KdV})$ a Crank-Nicholson scheme combined
with a relaxation method coming from Besse-Bruneau in
\cite{BB} and justified by Besse in \cite{Besse}. This type of scheme is of order two in space and time, which is appropriate for our purpose. 

\subsubsection{Numerical scheme for the KdV approximation}

\noindent Due to the identical structure of the two KdV equations of
$(\Sigma_{KdV})$, we only present the numerical
scheme for the first equation. Defining $u(t,x) = U_0(T,x-t)$, we can
reformulate this equation as follows
\begin{equation} \label{kdvnum}
\partial_t u + \partial_x u + \varepsilon \left[ \frac{3}{4} u
  \partial_x u + \frac{1}{6} \partial_x^3 u \right] = 0 \;.
\end{equation}
We use a Crank-Nicholson scheme 
and the relaxation method introduced by Besse-Bruneau in \cite{BB} and
justified by Besse in \cite{Besse} which replace the costly numerical
treatment of the nonlinear term by a predictive step. This provides us
with the following semi-discretized in time equation :
\begin{eqnarray*}
\frac{u^{n+1}-u^n}{dt} & + & \partial_x \left(\frac{u^{n+1}+u^n}{2}
\right) + \varepsilon \left[ \frac{3}{4} \left( \alpha u^{n+\frac{1}{2}}
    \partial_x \left(\frac{u^{n+1}+u^n}{2}\right) \right. \right. \\
  & + &  \left. \left. (1-\alpha) \frac{u^{n+1}+u^n}{2} \partial_x u^{n+\frac{1}{2}}
  \right) + \frac{1}{6} \partial_x^3 \frac{u^{n+1}+u^n}{2} \right] = 0\;,
\end{eqnarray*}
where the predictive term $u^{n+\frac{1}{2}}$ is defined as follows
$$
u^n = \frac{u^{n+\frac{1}{2}}+u^{n-1/2}}{2}\;.
$$
The discretization of the nonlinear term $u \partial_x u$ here takes
advantage of the two possible discretizations $u^{n+\frac{1}{2}}
    \partial_x \left(\frac{u^{n+1}+u^n}{2}\right)$ and
    $\frac{u^{n+1}+u^n}{2} \partial_x u^{n+\frac{1}{2}}$ by introducing a
    parameter $\alpha \in [0,1]$ and taking a convex combination of these
    possibilities.
Keeping in mind that we want to preserve the semi-discrete $L^2$ norm,
an easy integration by parts gives us the appropriate value $\alpha =
2/3$. We then choose the spatial discretization so that the
discrete $L^2$ norm is preserved by the scheme, which gives
the final discretization of (\ref{kdvnum}) :
\begin{equation}
\begin{array}{c}
\vspace{0.5em}
\frac{u_i^{n+1}-u_i^n}{\delta t} + \left( D_1
  \frac{u^{n+1}+u^{n}}{2}\right)_i + \varepsilon \left[ \frac{1}{4}
  \left(u_i^{n+\frac{1}{2}} + \frac{u_{i+1}^{n+\frac{1}{2}}+u_{i-1}^{n+\frac{1}{2}}}{2}\right)
\right. \\ 
\left. \left(D_1 \frac{u^{n+1}+u^{n}}{2}\right)_i 
 + \frac{1}{4} \frac{u_i^{n+1}+u_i^{n}}{2} \Big(D_1 u^{n+\frac{1}{2}} \Big)_i
 + \frac{1}{6} \left(D_3 \frac{u^{n+1}+u^{n}}{2}\right)_i \,\right] = 0\;,
\end{array}
\end{equation}
where the matrix $D_1$ and $D_3$ are to the classical centered
discretizations of the derivatives $\partial_x$ and $\partial_x^3$.\\
Once the solution of the KdV equation is obtained, the KdV approximation is built thanks to the relations (\ref{MB}), by computing the integrals corresponding the correcting terms with the composite trapezoidal quadrature rule. Here we have made the choice to save the solution $u$ of the KdV equation, and then build the approximate solution from it, which can be computationnaly costly in general, especially in 3-D. This is not problematic here since we work in 2-D. However, a optimised approach would be to see the equations (\ref{MB}) as an easy to discretize two dimensionnal pseudo-differential operator acting on the vector $(U_0,N_0)$. Multiplying the previous scheme with the discretized version of this operator would lead a new system to solve at each time step, which would give us directly the KdV approximation.

\vspace{1em}

\subsubsection{Numerical scheme for the Boussinesq system}

\noindent As far as the discretization of the Boussinesq system
$(\Sigma)$ is concerned, we consider the same ideas. Using a Crank-Nicholson
scheme and the same relaxation method, we aim here at preserving the
specific norm $|(v,\eta)|_{H^1_{\varepsilon}}^2 = |v|_{L^2}^2 + |\eta|_{L^2}^2 + \varepsilon a_2
|\partial_x v|_{L^2}^2 + \varepsilon a_4 |\partial
\eta|_{L^2}^2$. This quantity is indeed conserved by $(\Sigma)$ (see
\cite{BCL} for more details). To this end, the nonlinear terms $v
\partial_x v$, $\eta \partial_x \eta$, $ \eta \partial_x v$ and $v
\partial_x \eta$ are discretized in order to preserve both this
specific discrete norm and their symmetric structure.
Remarking that the equalities 
$$(v\partial_x v,v)_{L^2} = 0 \;\;;\;\; (\eta
\partial_x \eta, v)_{L^2} + (\eta \partial_x v, \eta)_{L^2} + (v
\partial_x \eta, \eta)_{L^2} = 0\;,
$$
hold for $(\Sigma)$ and using the
same kind of method as for the KdV equation leads to the following
semi-discretization of the nonlinear terms :
$$
\left\{
\begin{array}{l}
\vspace{0.5em}
v \partial_x v (n \delta t) \approx \frac{2}{3} v^{n+\frac{1}{2}} \partial_x
(\frac{v^{n+1}+v^{n}}{2}) + \frac{1}{3} \frac{v^{n+1}+v^{n}}{2}
\partial_x v^{n+\frac{1}{2}} \;, \\
\vspace{0.5em}
\eta \partial_x \eta (n \delta t) \approx \frac{2}{3} \eta^{n+\frac{1}{2}}
\partial_x (\frac{\eta^{n+1}+\eta^{n}}{2}) + \frac{1}{3}
\frac{\eta^{n+1}+\eta^{n}}{2} \partial_x \eta^{n+\frac{1}{2}} \;, \\
\vspace{0.5em}
\eta \partial_x v (n \delta t) \approx \frac{2}{3} \eta^{n+\frac{1}{2}}
\partial_x (\frac{v^{n+1}+v^{n}}{2}) + \frac{1}{3}
\frac{\eta^{n+1}+\eta^{n}}{2} \partial_x v^{n+\frac{1}{2}} \;, \\
v \partial_x \eta (n \delta t) \approx \frac{2}{3} v^{n+\frac{1}{2}} \partial_x (\frac{\eta^{n+1}+\eta^{n}}{2}) + \frac{1}{3} \frac{v^{n+1}+v^{n}}{2} \partial_x \eta^{n+\frac{1}{2}} \;.
\end{array}
\right.
$$ 
We then choose the spatial discretization so that the discrete
$H^1_{\varepsilon}$ norm is conserved, and these ruminations yield this final scheme :
\begin{equation}
\left\{
\begin{array}{c}
\vspace{0.5em}
\left((I-\varepsilon a_2 D_2) \frac{v^{n+1}-v^n}{\delta
    t}\right)_i + \left( (I-\frac{\varepsilon}{2} B) D_1
  \frac{\eta^{n+1}+\eta^{n}}{2}\right)_i \\
\vspace{0.5em} + \varepsilon \left[
  \left(M_1 \, \frac{v^{n+1}+v^{n}}{2}\right)_i  +  \left(M_2
    \, \frac{\eta^{n+1}+\eta^{n}}{2}\right)_i \right. \\ 
\vspace{1em} \left. +
  a_1 \left(D_3 \frac{\eta^{n+1}+\eta^{n}}{2}\right)_i \right] = 0 \\ 
\vspace{0.5em} \left((I-\varepsilon a_4 D_2) \frac{\eta^{n+1}-\eta^n}{\delta
    t}\right)_i + \left( (I-\frac{\varepsilon}{2} B) D_1
  \frac{v^{n+1}+v^{n}}{2}\right)_i \\ 
\vspace{0.5em} + \varepsilon \left[ \left(M_3
    \, \frac{\eta^{n+1}+\eta^{n}}{2}\right)_i  +  \left(M_4
    \, \frac{v^{n+1}+v^{n}}{2}\right)_i \right. \\ 
\left. + a_1
  \left(D_3 \frac{v^{n+1}+v^{n}}{2}\right)_i \right] = 0 \;,
\end{array}
\right.
\end{equation}
where the matrix $(M_i)_{1\le i \le 4}$ are as follows :
$$
\left\{
\begin{array}{l}
\vspace{0.5em}
\left(M_1 \, \frac{v^{n+1}+v^{n}}{2}\right)_i = \frac{1}{2} (v_i^{n+\frac{1}{2}} + \frac{v_{i+1}^{n+\frac{1}{2}}+v_{i-1}^{n+\frac{1}{2}}}{2}) \left(D_1 \frac{v^{n+1}+v^{n}}{2} \right)_i 
\\ \vspace{0.5em}  \hspace{11em} + \frac{1}{2} 
\left(D_1 v^{n+\frac{1}{2}}\right)_i \frac{v_i^{n+1}+v_i^{n}}{2} \;, \\ \vspace{0.5em}
\left(M_2 \, \frac{\eta^{n+1}+\eta^{n}}{2}\right)_i = \frac{1}{3} \eta_i^{n+\frac{1}{2}} \left(D_1 \frac{\eta^{n+1}+\eta^{n}}{2} \right)_i 
\\ \vspace{0.5em} \hspace{11em} + \frac{1}{6} 
\left(D_1 \eta^{n+\frac{1}{2}}\right)_i \frac{\eta_i^{n+1}+\eta_i^{n}}{2} \;, \\ \vspace{0.5em}
\left(M_3 \, \frac{v^{n+1}+v^{n}}{2}\right)_i = \frac{1}{3} \frac{\eta_{i+1}^{n+\frac{1}{2}}+\eta_{i-1}^{n+\frac{1}{2}}}{2} \left(D_1 \frac{v^{n+1}+v^{n}}{2} \right)_i 
\\ \vspace{0.5em}  \hspace{11em} + \frac{1}{3} 
\left(D_1 \eta^{n+\frac{1}{2}}\right)_i
\frac{\eta_{i+1}^{n}+\eta_{i-1}^{n}-\eta_i^n}{2} \;, \\ \vspace{0.5em}
\left(M_4 \, \frac{\eta^{n+1}+\eta^{n}}{2}\right)_i = \frac{1}{6} (v_i^{n+\frac{1}{2}} + \frac{v_{i+1}^{n+\frac{1}{2}}+v_{i-1}^{n+\frac{1}{2}}}{2}) \left(D_1 \frac{v^{n+1}+v^{n}}{2} \right)_i 
\\ \vspace{0.5em} \hspace{11em} + \frac{1}{6} 
\left(D_1 v^{n+\frac{1}{2}}\right)_i \frac{v_i^{n+1}+v_i^{n}}{2} \;. \\
\end{array}
\right.
$$
The matrix $D_1$ and $D_3$ are as defined in the KdV scheme, and the matrix $D_2$ is the classical centered discretization of the derivative $\partial_x^2$.

\vspace{1em}

\subsubsection{Initial data}
 
\noindent Let us now talk about the initialization of the two schemes. First, all the prevision terms are initialized with a simple explicit integration of the equations on a half-step in time. Then, the initial conditions are chosen such that the simulated wave is unidirectional and propagating to the right. To this end, we first take the initial data of the second KdV equation to be zero. Then the system  $(\Sigma_{KdV})$ reduces to the equation (\ref{kdvnum}) for which we know the existence of solitary waves expressed as follows :
\begin{equation} \label{soluce}
u(t,x) = \frac{\alpha}{\cosh^2\big(k(x-c t+l)\big)}\;,
\end{equation}
with $c=1+\frac{\varepsilon \alpha}{4}$, $k=\sqrt{\frac{3\alpha}{8}}$ and $\alpha$, $l$ being arbitrary.\\
It is hence natural to specify the initial condition for the KdV equation (\ref{kdvnum}) as follows :
\begin{equation} \label{init}
u(t=0,x) = u_0(x) = \frac{\alpha}{\cosh^2\big(k(x+d)\big)}\;.
\end{equation}
Finally, and because of the way the KdV approximation was constructed from the Boussinesq model, we specify the initial conditions for this latter as follows :
$$
v(t=0,x) = \eta(t=0,x) = \frac{1}{2}\, u_0(x) \;.
$$

\subsubsection{Validation of the numerical method}

With the initial data  (\ref{init}), the KdV scheme is expected to
propagate the corresponding solitary wave to the right, without any
deformation for all time.  In order to validate this scheme,
the numerical results obtained with the
initial data (\ref{init}) have been compared with the analytical solution
(\ref{soluce}). The following relative
errors on the free surface have been computed in the $L^{\infty}$ norm for several values
of epsilon and for computation times $T = 1/\varepsilon$ :

\begin{table}[h]
\vspace*{1em}
\begin{center}
\begin{tabular}{|*{7}{c|}}
\hline
$\varepsilon$ & $T$ & $L$ & $\delta x$ & $\delta t$ & relative error \\
\hline
0.05 & 20 & 80 & 0.03 & 0.03 & {\small $1.5546.10^{-3}$} \\
\hline
\vspace*{0.05em}
0.1 & 10 & 80 & 0.04 & 0.04 & {\small $1.3717.10^{-3}$} \\
\hline
\vspace*{0.05em}
0.2 & 5 & 80 & 0.05 & 0.05 &  {\small $1.0534.10^{-3}$} \\
\hline
\end{tabular}
\end{center}
\end{table}

\noindent where $L$ is the length of the computational domain and $\delta
x$,$\delta t$ are respectively the spatial and time discretization
steps. These results allow to validate the scheme proposed for the KdV equations.

\subsection{Numerical results and comments}

\subsubsection{Numerical results}

As specified in \cite{FC}, the choice of the parameters $a_1,a_2,a_4$ is very interesting in a numerical point of view. Indeed, the parameter $a_1$ controls the presence of the dispersive terms $\partial_x^3 v$ and $\partial_x^3 \eta$ whereas the parameters $a_2$ and $a_4$ correspond to the terms $\partial_x^2\partial_t v$ and $\partial_x^2\partial_t \eta$. These last terms have the main advantage of being regularizing terms analytically and numerically speaking, they smooth in some way the solution because they provide a control of the quantities $\partial_x v$ and $\partial_x \eta$ in the $L^2$ norm. We decided to use here the system $(\Sigma)$ corresponding to $a_1=a_2=a_4=1/12$ because it is likely to provide the better results.

\noindent All the forthcoming results are
expressed in non-dimensionalized variables.  We recall that both the
free surface and the bottom are of size $\varepsilon$ :
$y=\varepsilon\eta$ for the free surface and $y=-1+\varepsilon b$
for the bottom. However, in order to get clear and readable results,
we have plotted a rescaled free surface $y=\eta$ and a rescale bottom
$y=-1 + b$. A quick word on
the duration $T$ of the simulations : the 
previous secion and \cite{FC} provide us with a justification of
the models on large time scales of order $O(1/\varepsilon)$. We have
decided - only in the first example of the step - to overtake this large
time scale and simulate the models on the very large time
$T=1/\varepsilon^{3/2}$, in order  
to see if the model remains stable on such time scales. 

\noindent The three models have been tested on two different examples of
bottom. The first one correspond to a step at the bottom, defined
similarly to \cite{Dutykh} by 
\begin{equation} \label{stepn}
b(x) = \;
\left\{
\begin{array}{l}
\vspace{0.25em}
0\;,\;\;\forall x \in \left[0,\frac{L}{2}-\frac{3}{2}\right] \;,\\
\vspace{0.25em}
\frac{\beta_0}{2}\left(1+\sin\Big(\frac{\pi}{3}(x-\frac{L}{2})\Big)\right) \;,\;\;\forall x \in \left[\frac{L}{2}-\frac{3}{2},\frac{L}{2}+\frac{3}{2}\right] \;,\\
\beta_0\;,\;\;\forall x \in \left[\frac{L}{2}+\frac{3}{2},L\right] \;,
\end{array}
\right.
\end{equation}
where $\beta_0$ is an arbitrary constant of order $O(1)$ and $L$ is the length of the computation domain.
The second example corresponds to a slowly varying sinusoidal bottom, defined as :
\begin{equation} \label{slow}
b(x) = b_0\sin\left(\frac{\pi}{2}+\frac{2\pi}{l}
  x\right)\;,\;\;\forall x \in \R\;,
\end{equation}
where $l$ is defined by
$l=\frac{1+\varepsilon\alpha/4}{\varepsilon}$ and $\alpha$
is the amplitude of the initial data defined in (\ref{init}). 

\vspace{1em}

\noindent The following results show the snapshots of the simulations
at different times - so that the time evolution
is relatively visible - and the evolution of the relative $L^{\infty}$ error between
the free surfaces obtained with the Boussinesq model and respectively the
KdV approximation and the topographically modified approximation. 
The three models have been systematically plotted together in the same pictures in
order to compare efficiently their respective behaviours. The numerical simulations have been performed for
different values of $\varepsilon$ in the case (\ref{stepn}) of a step : $\varepsilon = 0.05$ and
$\varepsilon = 0.2$, which are typical values of the upper part of the range of validity of
the long waves approximation. As far as the case (\ref{slow}) is
concerned, we simulated the models for the value $\varepsilon=0.1$. For all the simulations, the amplitude $\alpha$ of the initial free surface and the constant $\beta_0$ linked to the bottom have been taken equal to $0.5$. Here is a global tabular precising all the values of interest used in the simulations.

\begin{table}[ht]
\vspace*{1em}
\begin{center}
\begin{tabular}{|*{7}{c|}}
\hline
Figure & Bottom & $\varepsilon$ & $T$ & $L$ & $\delta x$ & $\delta t$ \\
\hline
\ref{figure1} & step & 0.05 & 89 & 140 & 0.03 & 0.03 \\
\hline
\ref{figure5} & step & 0.2 & 12 & 80 & 0.05 & 0.05 \\
\hline
\ref{figure9} & sinusoidal & 0.1 & 10 & 20 & 0.04 & 0.04 \\
\hline
\end{tabular}
\end{center}
\end{table}

\noindent The figures
\ref{figure2}, \ref{figure6}, \ref{figure10}
show the relative error between the computed free surfaces of the
different models for each value of $\varepsilon$ and for the two cases of bottom.

\begin{figure}
  \caption{Influence of the step for $\varepsilon=0.05$}
  \label{figure1}
  \centering
  \epsfig{file=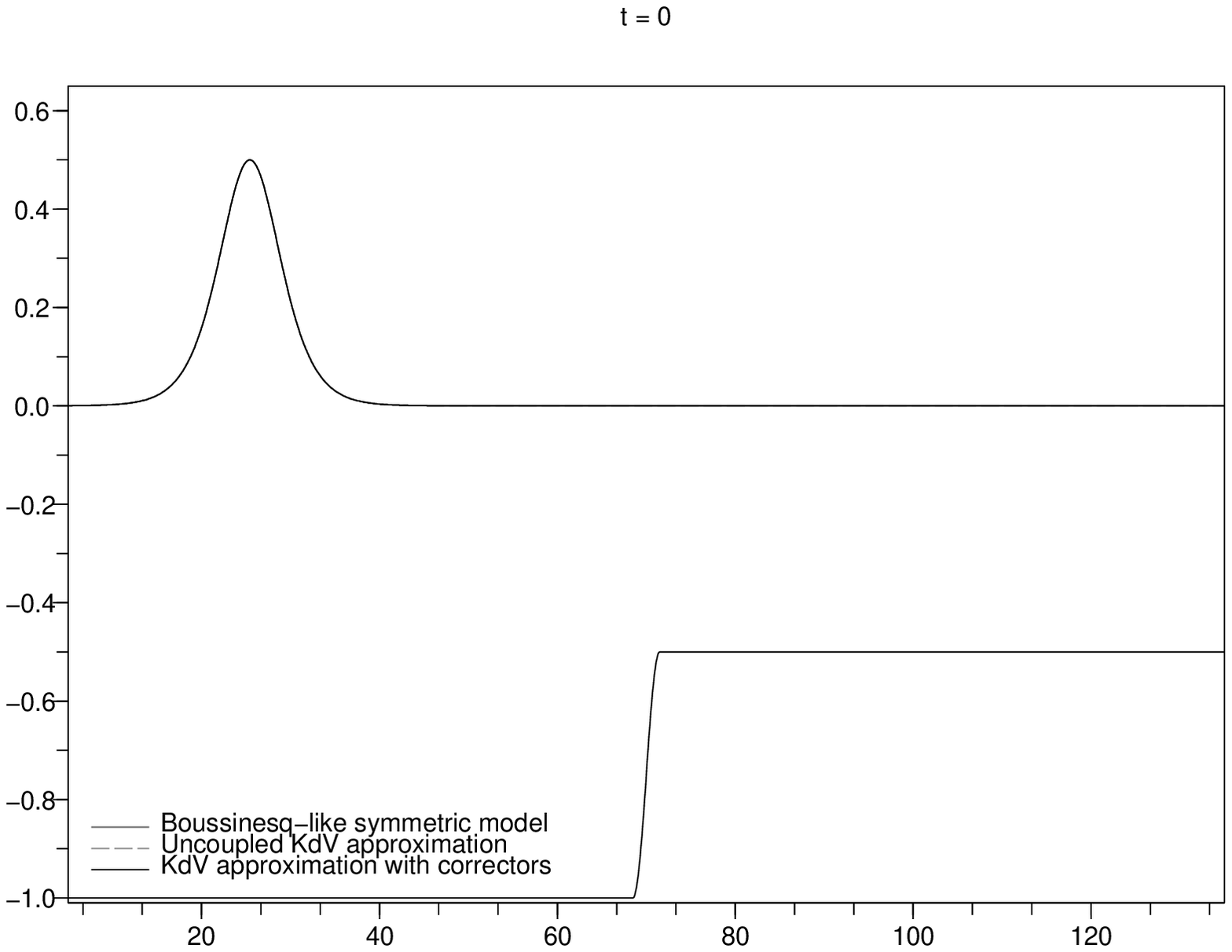,width=384pt,height=272pt}
  \epsfig{file=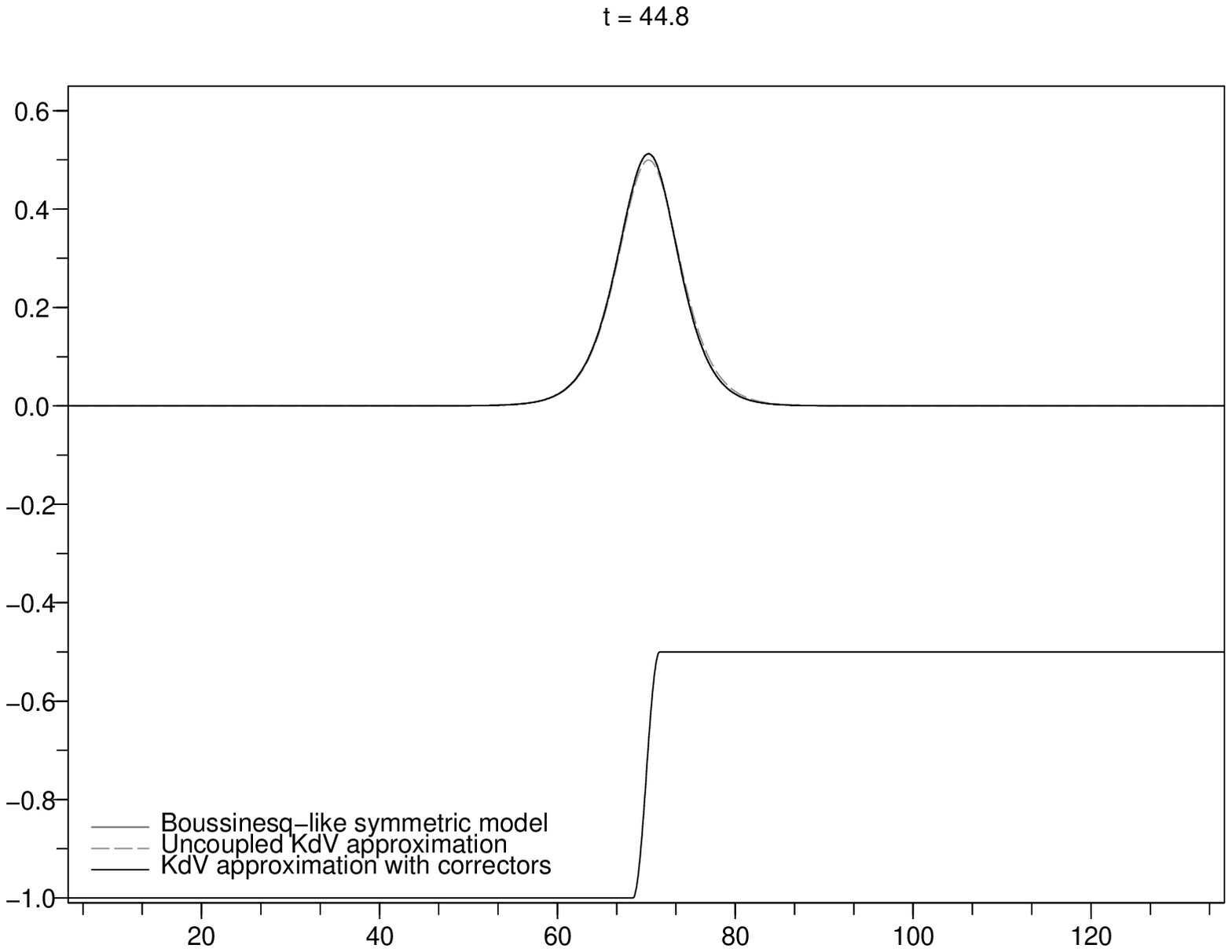,width=384pt,height=272pt}
\end{figure}
\begin{figure}
  \centering
  \epsfig{file=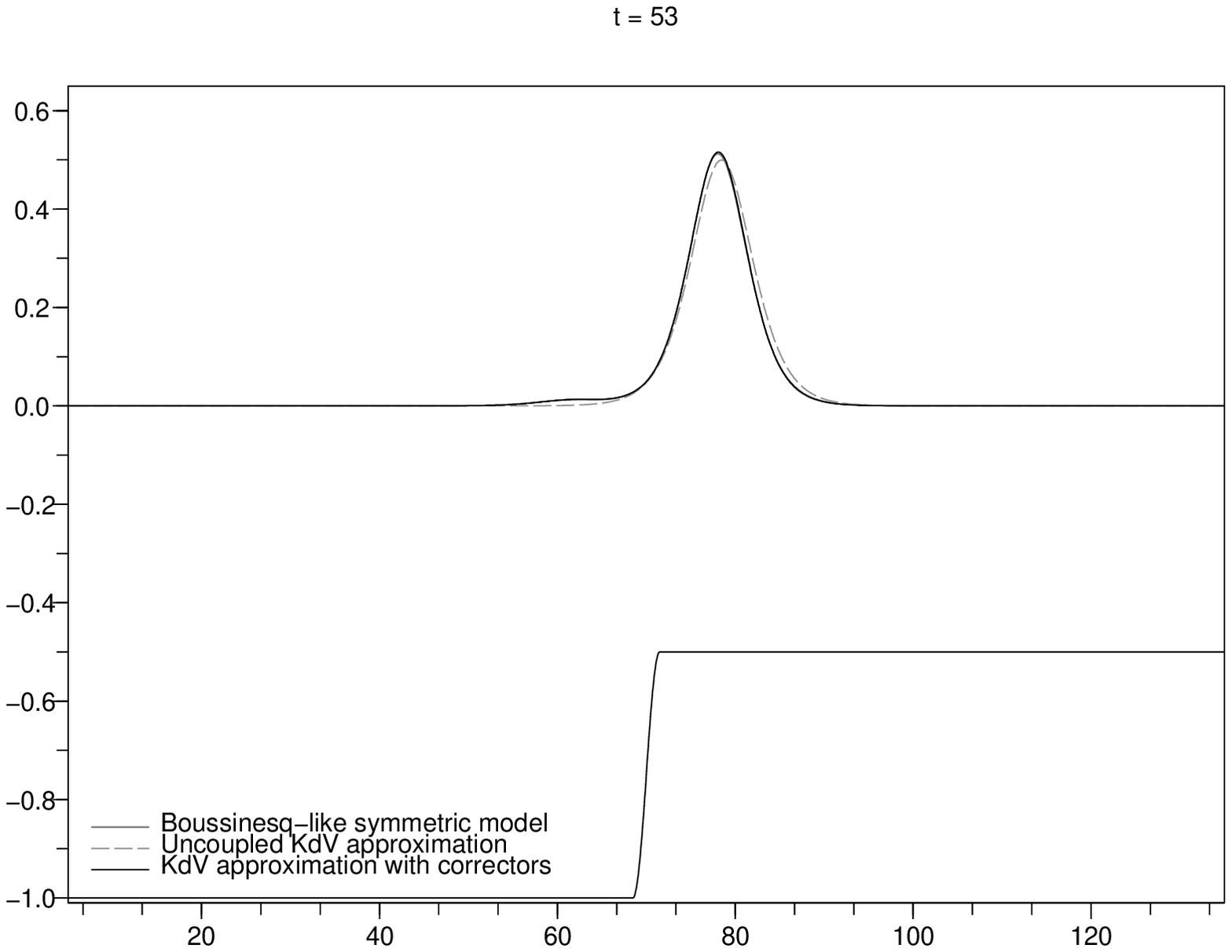,width=384pt,height=280pt}
  \epsfig{file=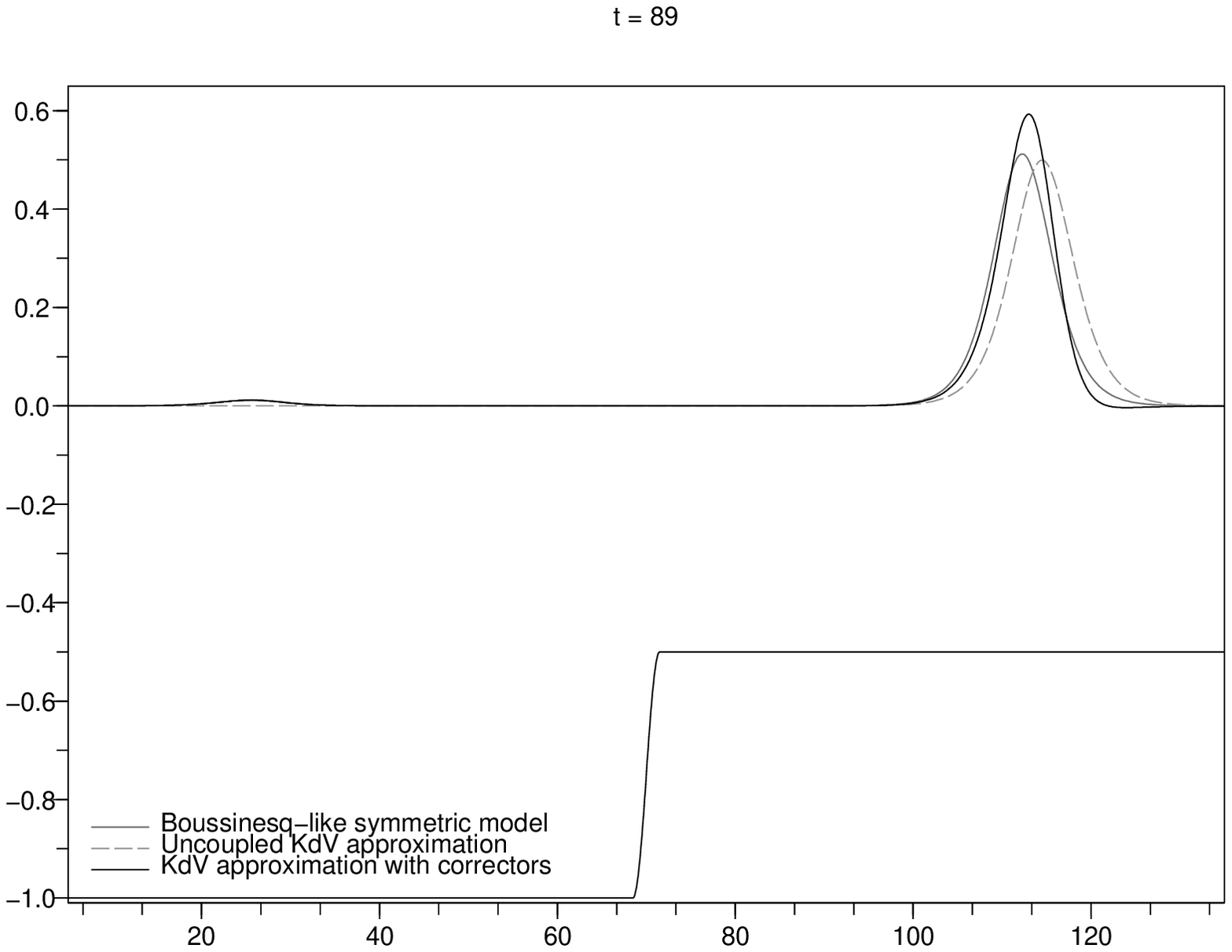,width=384pt,height=280pt}
\end{figure}
\begin{figure}
  \caption{Relative $L^{\infty}$ error between the free surfaces
    for $\varepsilon=0.05$}
  \label{figure2}
  \centering
  \epsfig{file=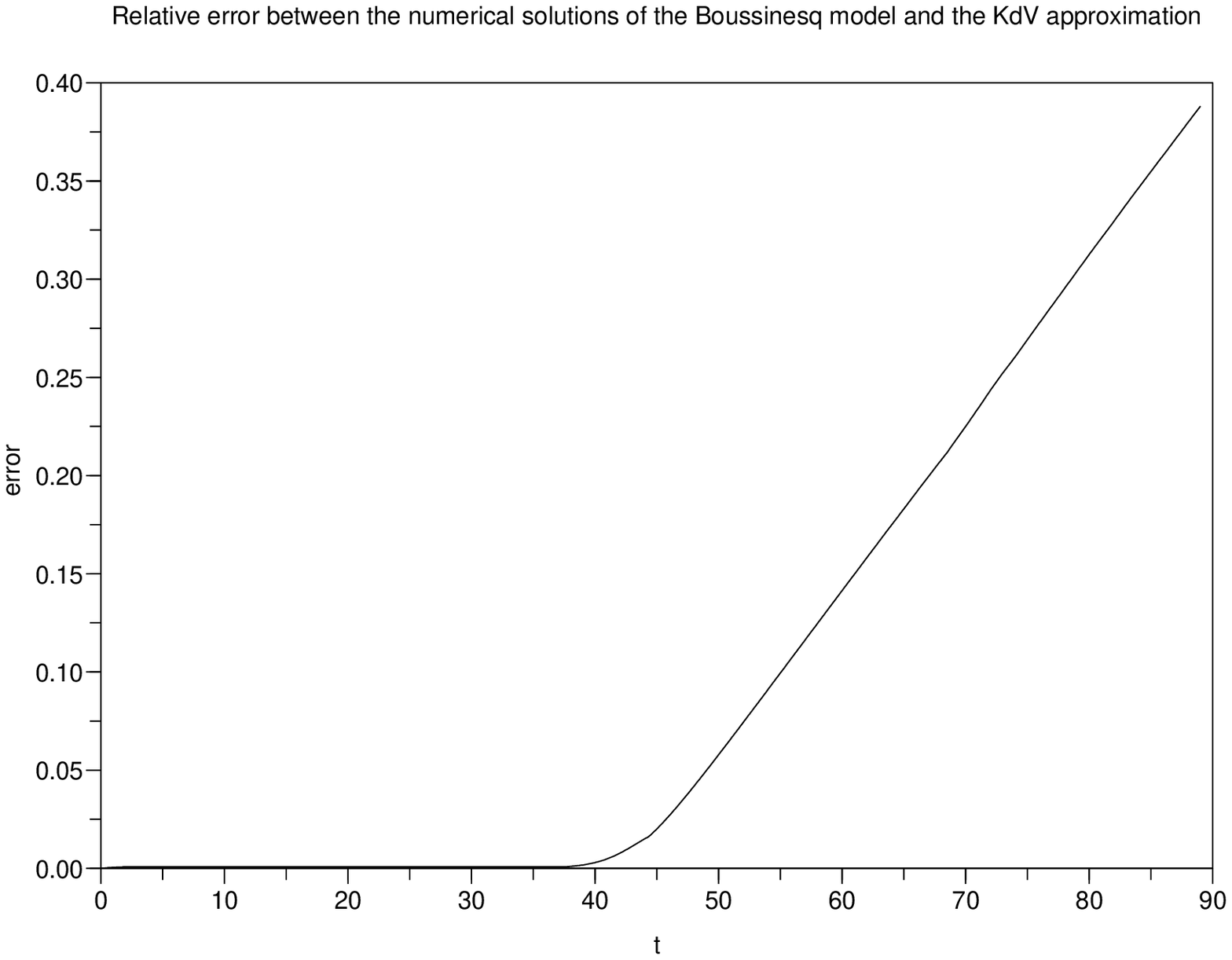,width=384pt,height=272pt}
  \epsfig{file=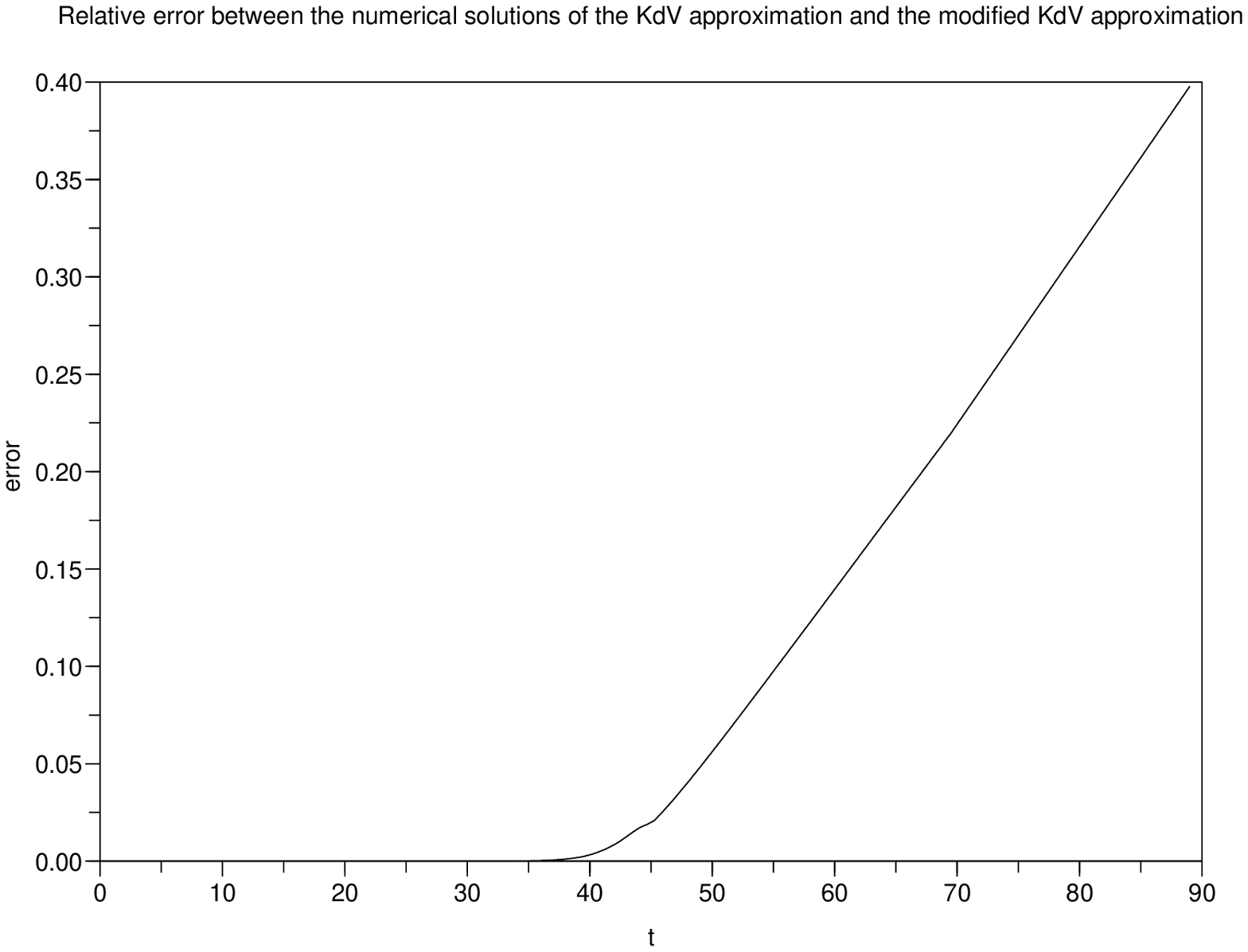,width=384pt,height=272pt}
\end{figure}
\begin{figure}
  \caption{Influence of the step for $\varepsilon=0.2$}
  \label{figure5}
  \centering
  \epsfig{file=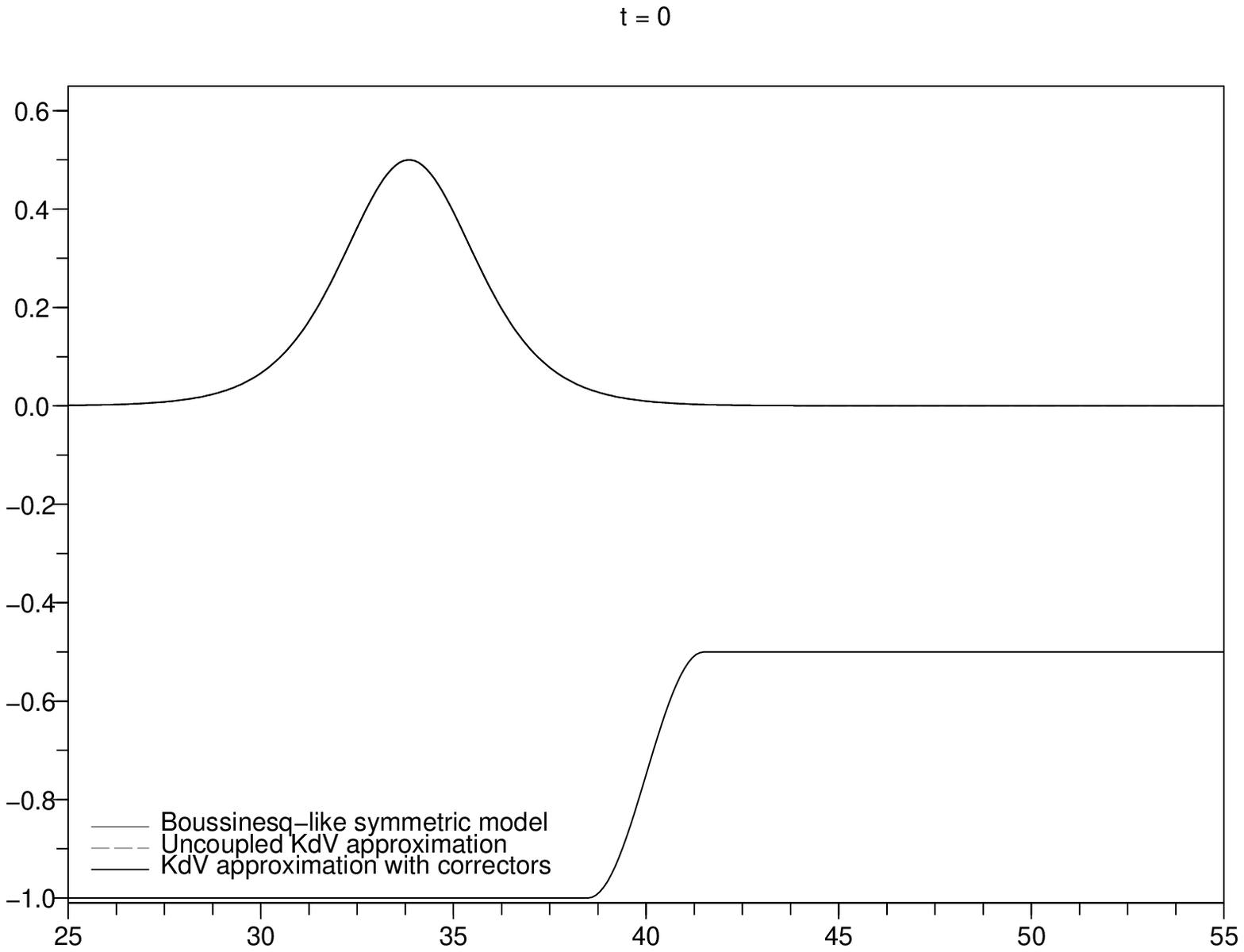,width=384pt,height=272pt}
  \epsfig{file=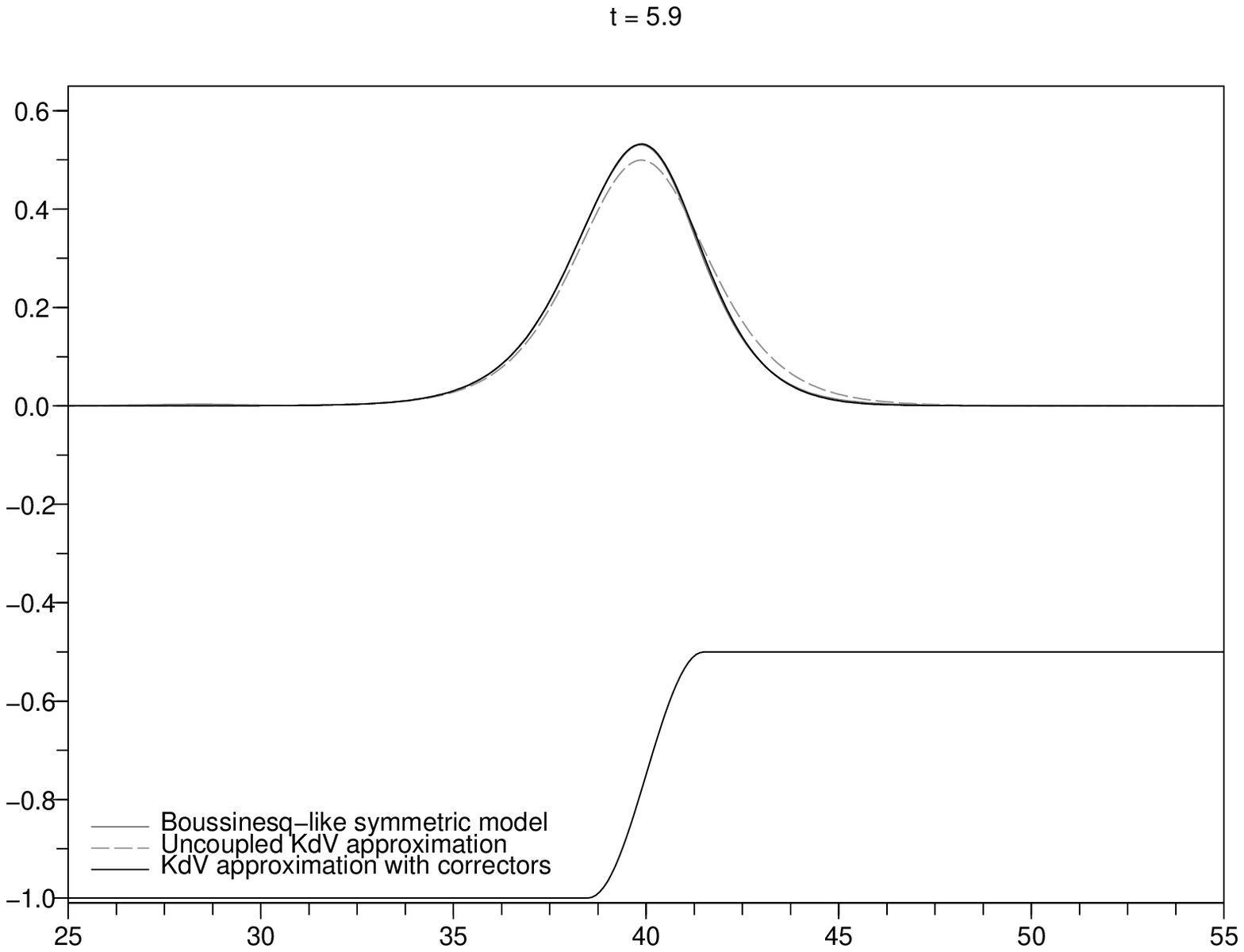,width=384pt,height=272pt}
\end{figure}
\begin{figure}
  \centering
  \epsfig{file=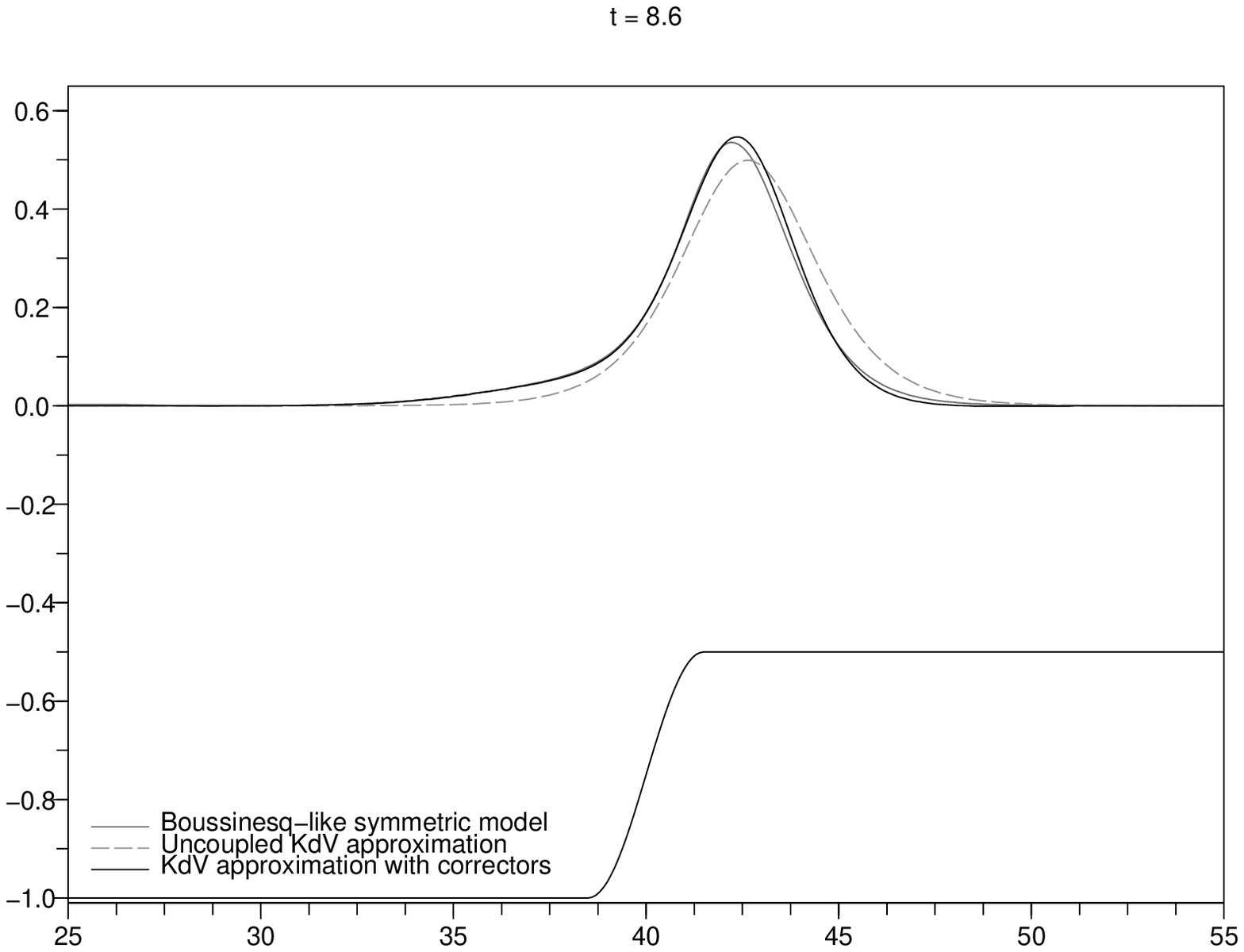,width=384pt,height=280pt}
  \epsfig{file=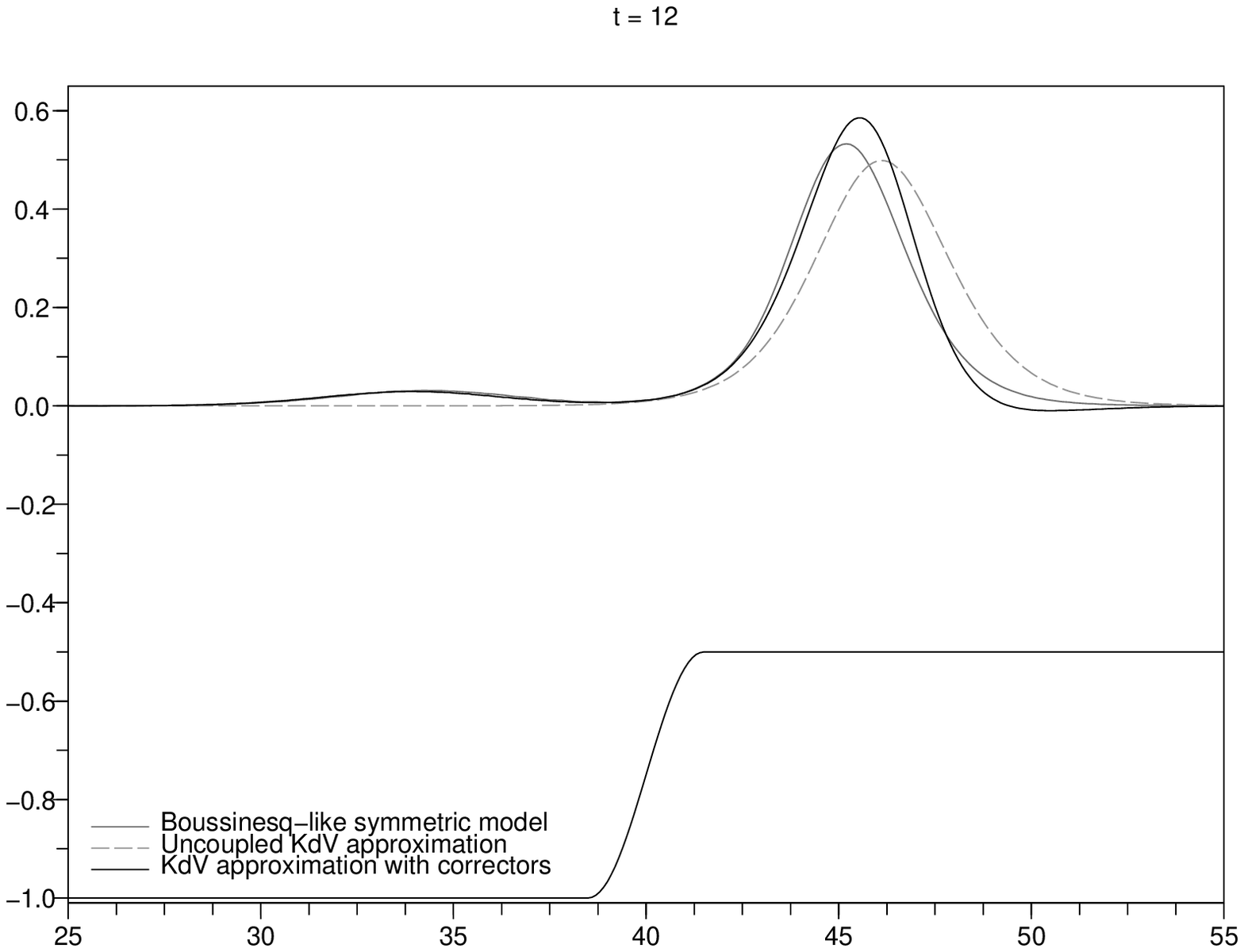,width=384pt,height=280pt}
\end{figure}
\begin{figure}
  \caption{Relative $L^{\infty}$ error between the free surfaces
    for $\varepsilon=0.2$}
  \label{figure6}
  \centering
  \epsfig{file=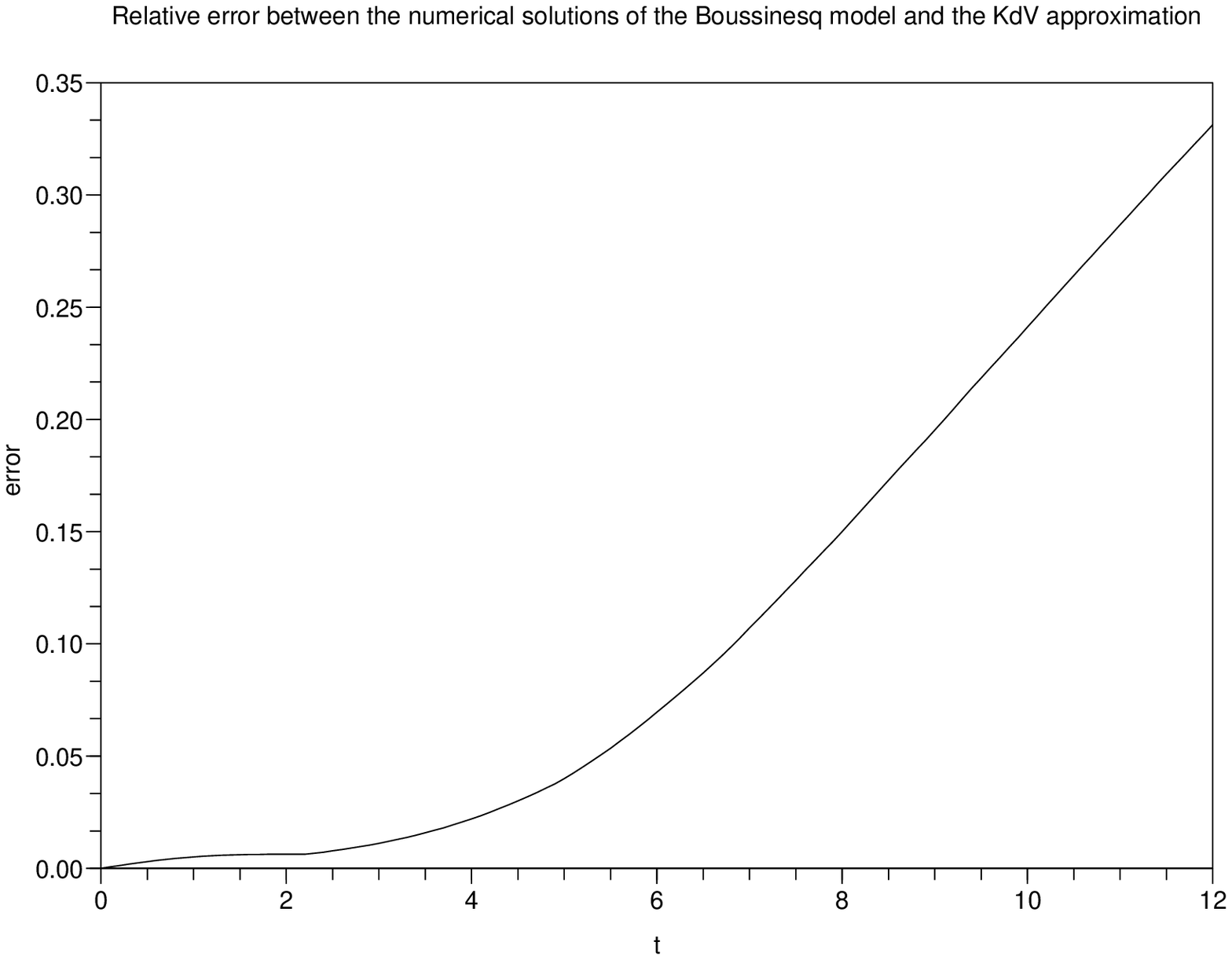,width=384pt,height=272pt}
  \epsfig{file=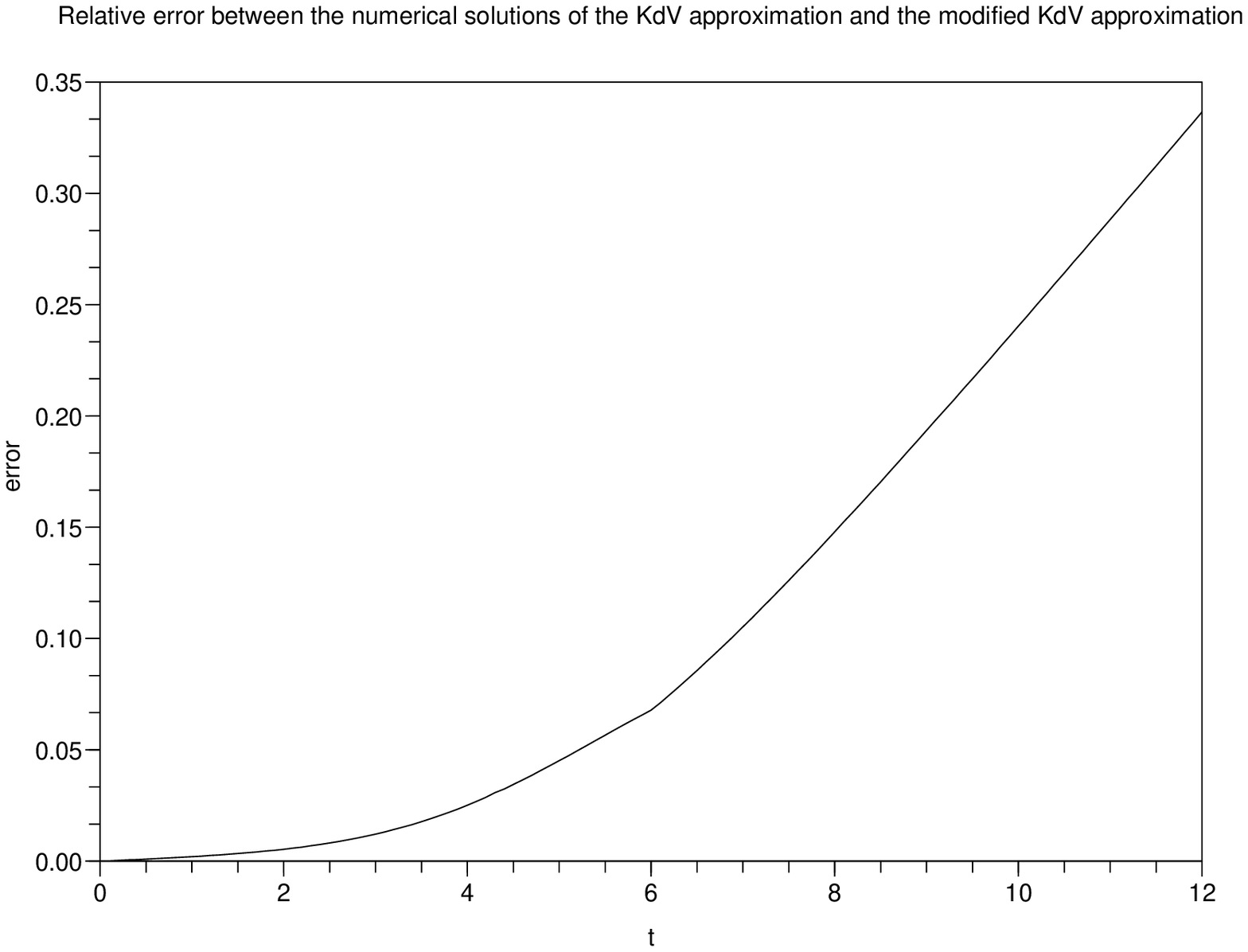,width=384pt,height=272pt}
\end{figure}
\begin{figure}
  \centering
  \caption{Influence of a slow sinusoidal bottom for
    $\varepsilon=0.1$}
  \label{figure9}
  \epsfig{file=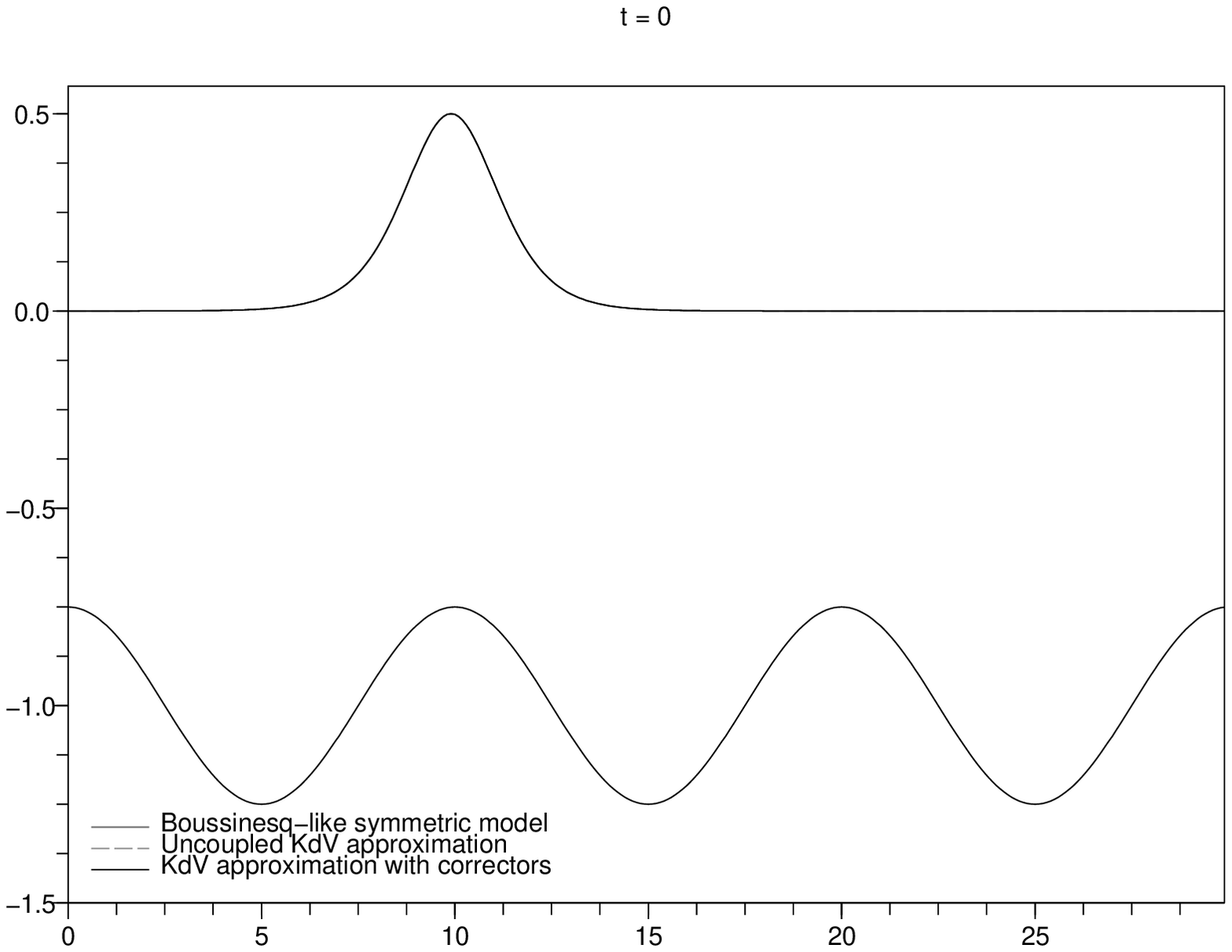,width=384pt,height=270pt}
  \epsfig{file=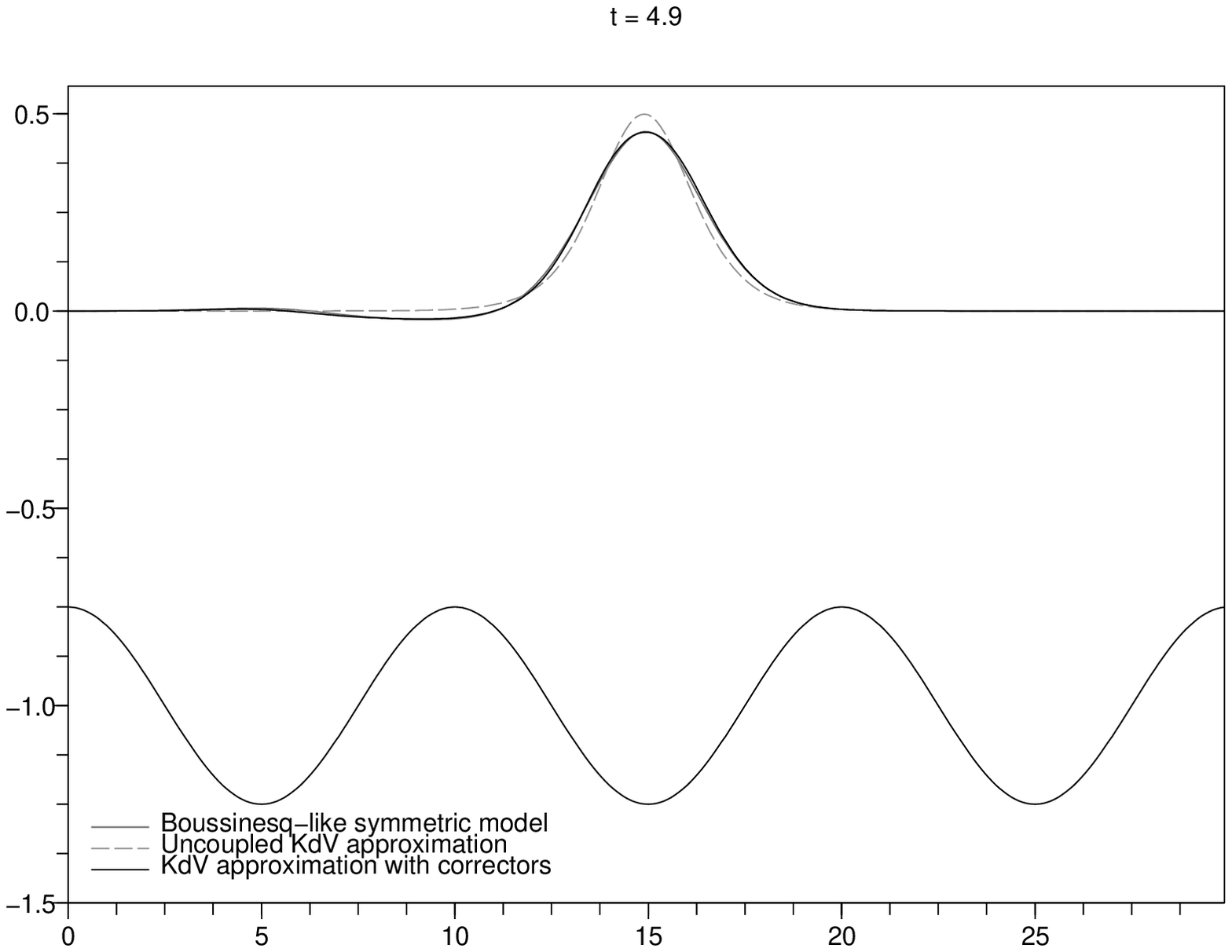,width=384pt,height=270pt}
\end{figure}
\begin{figure}
  \centering  
  \epsfig{file=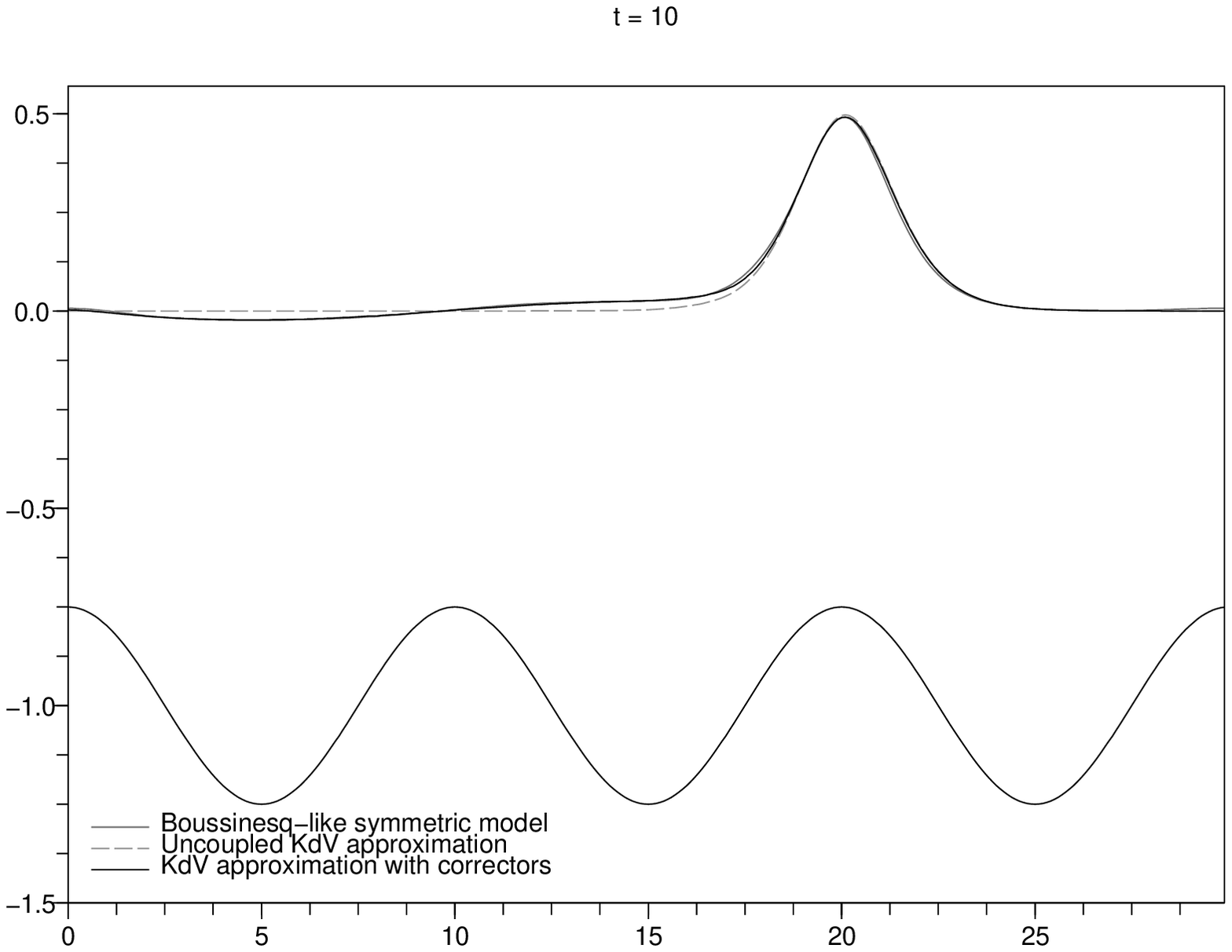,width=384pt,height=270pt}
  \caption{Relative $L^{\infty}$ error between the free surfaces
    for $\varepsilon=0.1$}
  \label{figure10}
  \centering
  \epsfig{file=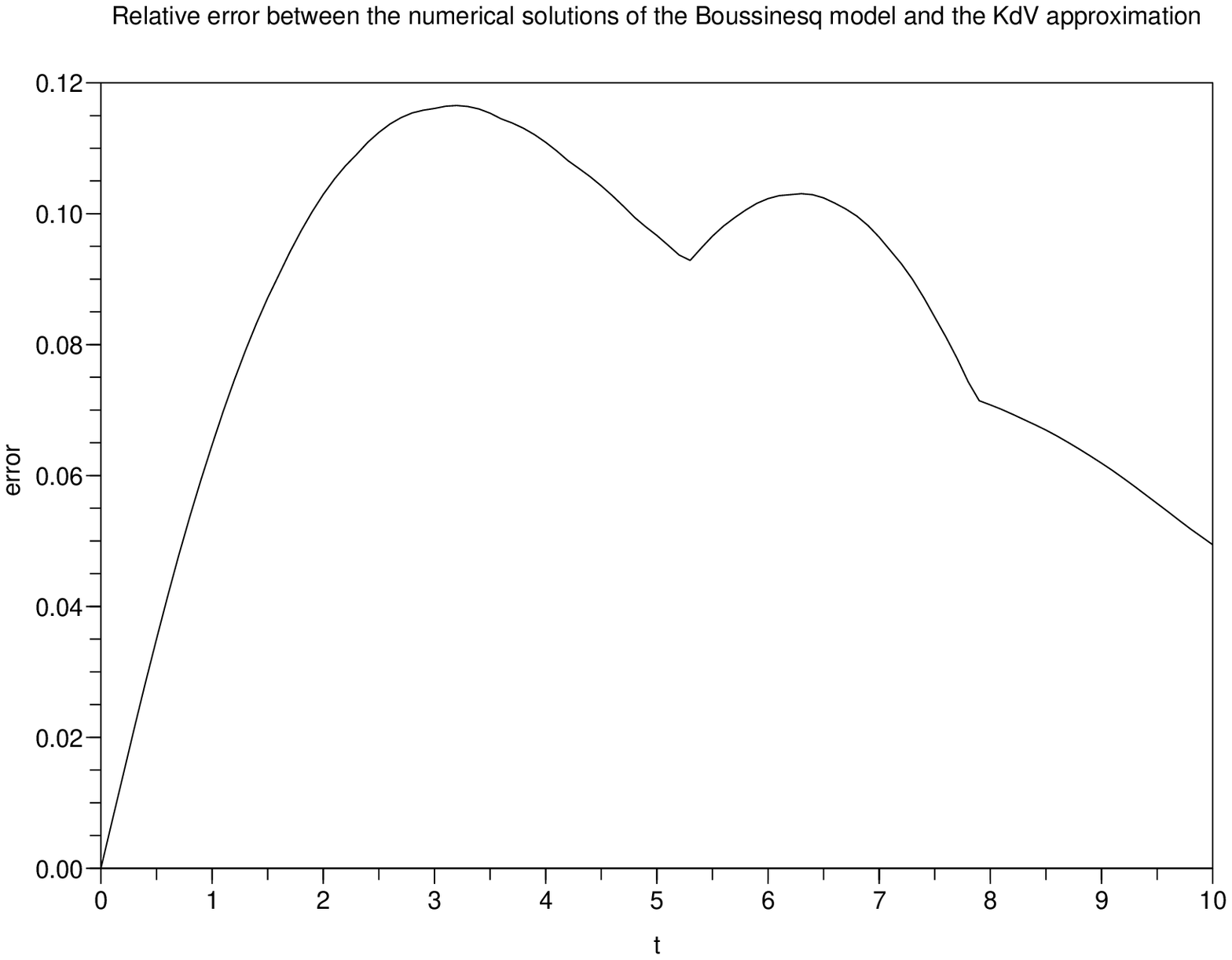,width=384pt,height=270pt}
\end{figure}
\begin{figure}
  \epsfig{file=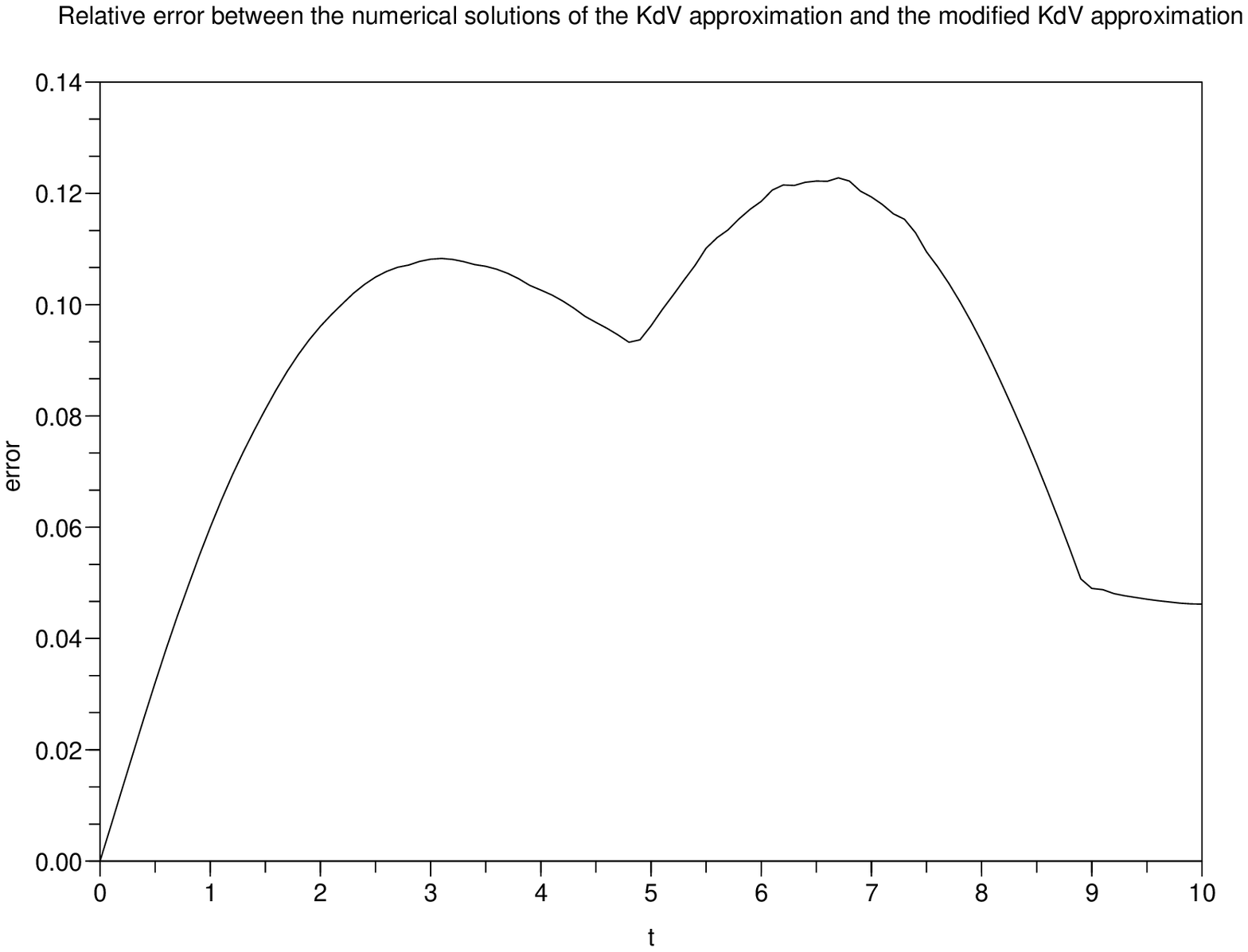,width=384pt,height=270pt}
\end{figure}

\vspace*{2em}

\subsubsection{Comments}

For the sake of readability, we call the computed waves as follows :
$\mathcal{B}$ denotes the wave coming from the Boussinesq model
$(\Sigma)$, $\mathcal{K}_{topo}$ denotes the wave produced by the
topographically modified KdV approximation $(\mathcal{M}_b)$, and
finally $\mathcal{K}$ denotes
the solitary wave resulting from the usual KdV approximation
$(\mathcal{M})$. 

\vspace{1em}

In the case of the step, we observe in figures \ref{figure1} and \ref{figure5} 
that for all tested values of
$\varepsilon$, both the Boussinesq model and alternative version
of the KdV approximation succeed in reproducing the phenomenon of
reflexion on the bathymetry : a smaller solitary wave appears on the third snapshots when the
main wave goes over the step. This reflected wave propagates to the
left at the same speed as the main wave. The classical uncoupled KdV
model cannot - of course - reproduce this phenomenon seing that it
does not depend at all on the bottom topography. Moreover, the
Boussinesq and topographically modified KdV models successfully
describe the following expected physical phenomenons : the shoaling
which corresponds to the growth in amplitude of the wave after the
step ; the deceleration of the wave after the step : the waves
$\mathcal{B}$ and $\mathcal{K}_{topo}$ are behind the wave
$\mathcal{K}$ which propagates at a constant speed on the last snapshots of 
figures \ref{figure1} and \ref{figure5} ; and
finally the loss of symmetry and the narrowing of the wave, which can
be remarked by comparing - f.e. on the last snapshot of figure \ref{figure5} -
the distances between several points of the
resulting waves at different heights : the solitary wave $\mathcal{K}$
propagates without any deformation and
remains symmetric, whereas the symmetry and width of $\mathcal{B}$ and
$\mathcal{K}_{topo}$ are modified by the step. All these phenomenons are the premisses of the
process of wave breaking, and they are all successfully reproduced by
the Boussinesq and new KdV models.

\vspace{1em}

\noindent A very interesting remark on our KdV model is that the role
of each correcting term in the approximation can be intuitively
identified, and these intuitions have been confirmed with several
simulations - that are not presented here - in the case of the
step. Indeed, it appears clearly on figures \ref{figure1} and \ref{figure5} that :
\begin{itemize}
\item the correcting term $\frac{\varepsilon}{8} \int_0^t \partial_x b
(x-t+s) U_0(x-t+2s) ds$ is responsible for the birth and propagation
of the reflected wave, and for the very beginning of the shoaling effect,
\item the correcting term $\frac{\varepsilon}{8} U_0(x-t)
  (b(x)-b(x-t))$ clearly reproduces the pursuit of the shoaling after
  the step,
\item the correcting term $\frac{\varepsilon}{4} \partial_x U_0(x-t)
  \int_0^t b (x-t+s) ds$ is reponsible for the deceleration and loss of width and symmetry  
of the main wave after the step. 
\end{itemize}

\vspace{1em}

\noindent In the case of a slowly varying sinusoidal bottom - corresponding to figure \ref{figure9} -
the effect of the bottom is also clearly visible on $\mathcal{B}$ and
$\mathcal{K}_{topo}$. To understand the speed variations of these waves, we
have to keep in mind that the comparison is made with the wave
$\mathcal{K}$ which evolves as if the
bottom was flat and located at the height $y=-1$. Consequently, 
when $\mathcal{B}$ and
$\mathcal{K}_{topo}$ propagate above the downward part of the
sinusoidal gap - see second snapshot of \ref{figure9} - they cross two different areas : a first area - for a
time $t \le T/4$ - where
the depth is lower than for a flat bottom located at $y=-1$, and a
second area - for $t \le T/2$ - where this is the contrary. This
explains why the waves $\mathcal{B}$ and
$\mathcal{K}_{topo}$ are located for $t=T/2$ at the same
position as $\mathcal{K}$ : the waves $\mathcal{B}$ and
$\mathcal{K}_{topo}$ have
speeded up over the first area in comparison with $(\mathcal{K})$, and
then have decelerated over the second area. However, we can see - still on the second snapshot of \ref{figure9} - 
that these waves are larger than $(\mathcal{K})$ at $t=T/2$ : this is due
to the loss of amplitude of the waves $\mathcal{B}$ and
$\mathcal{K}_{topo}$ during this downward part which makes the waves
be naturally wider. About this loss of amplitude, it is explained by
the fact that the amplitude of the bottom decreases over the downward
part of the sinusoid, which produces the inverse effect of the
shoaling from the previous case of the step. In addition, we can see a reflected
wave for $\mathcal{B}$ and $\mathcal{K}_{topo}$ - a depression this
time - which goes to the left at the same 
speed as the main wave, which corresponds to the same phenomenon of
bathymetric reflexion as in the case of the step : the amplitude of
the bottom decreases and thus produces a depression wave that
propagates in the opposite direction. We can see that this depression
wave is here larger than for a step because of the slow variations of
the bottom. As far as the upward part - see last snapshot of \ref{figure9} - of the sinusoid is concerned, all the
previously described effects happen in an inverted way, and we
finally recover three identical main waves for the three models. At
this final point, the only remaining visible effects of the crossed topography
are the reflected waves.

\vspace{1em}

\noindent As specified earlier, we decided to simulate the models on
the very large time $T=1/\varepsilon^{3/2}$ in the case of the
step on figures \ref{figure1} and \ref{figure5}. All these models have been proved to be valid on the time scale
$O(1/\varepsilon)$ and it is interesting to check numerically their
validity - or not - on larger time scales. 
For a time $T=1/\varepsilon^{3/2}$ , we can observe on the last two snapshots of figures \ref{figure1} and \ref{figure5} that the wave $\mathcal{K}_{topo}$ goes on growing more and more in amplitude, and
that a depression wave deepens in front of the main wave. These
effects are obviously not physical and can be explained by the fact
that the size of the correcting term $\frac{\varepsilon}{4} \partial_x U_0(x-t)
  \int_0^t b (x-t+s) ds$ evolves in time like $\varepsilon/t$ as we
  saw in the previous section : on a time $T=1/\varepsilon^{3/2}$,
  this size become of order $O(1/\sqrt{\varepsilon})$, which explains
  why this model diverges from the other models on this time
  scale. This is the main restriction of this model, in comparison
  with the Boussinesq one which seems to remain stable on very large
  time scales. An interesting perspective would be to look for higher order terms 
  - like Wright in \cite{Wright} - in the approximation to deal with this problem. 

\vspace{1em}

\noindent To sum up, the results on these two examples of bottom show that both Boussinesq and
topographically modified KdV models are able to reproduce the expected
physical phenomenons : reflection, shoaling, loss of speed and
symmetry. This is of course not the case for the usual KdV
approximation which is independent from the bottom topography. Even
if we can isolate the role of each correcting terms with our modified KdV
aproximation, this one diverges when time goes over the theoretical
limit time of validity $T=1/\varepsilon$. The Boussinesq does not have this drawback and remains
stable on time scales of order $O(1/\varepsilon^{3/2})$.

\vspace{3em}

\noindent {\bf Acknowledgments}. 
{\it This work was supported by the ACI Jeunes chercheurs du minist\`ere de
la Recherche ``Dispersion et nonlin\'earit\'es''.}

\end{document}